\documentclass[12pt]{article}

\usepackage{latexsym,amsfonts,amssymb,exscale,enumerate}
\usepackage{amsmath}
\usepackage{amsthm}

\usepackage{graphicx}
\usepackage[all]{xy}
\SelectTips{cm}{}

\newcommand{\seq}{{\rm Seq}}
\renewcommand{\bot}{{\rm bot}}
\renewcommand{\top}{{\rm top}}
\newcommand{\Pol}{\cal{P}o\ell}
\newcommand{\ii}{ \textbf{\textit{i}}}
\newcommand{\jj}{ \textbf{\textit{j}}}
\newcommand{\kk}{ \textbf{\textit{k}}}
\newcommand{\jRi}{{_{\jj}R(\nu)_{\ii}}}
\newcommand{\iRi}{{_{\ii}R(\nu)_{\ii}}}
\newcommand{\iRj}{{_{\ii}R(\nu)_{\jj}}}
\newcommand{\jBi}{{_{\jj}B_{\ii}}}

\newcommand{\jSi}{{_{\jj}S_{\ii}}}
\newcommand{\joi}{{_{\jj}1_{\ii}}}

\newcommand{\BNC}{NH}

\newcommand{\refequal}[1]{\xy {\ar@{=}^{#1}
(-1,0)*{};(1,0)*{}};
\endxy}

\newcommand{\U}{\dot{{\rm U}}}
\def\sbinom#1#2{\left( \hspace{-0.06in}\begin{array}{c} #1 \\ #2
\end{array}\hspace{-0.06in} \right)}

\usepackage{fancyheadings}
\newcommand{\Hom}{{\rm Hom}}
\newcommand{\HOM}{{\rm HOM}}
\renewcommand{\to}{\rightarrow}
\newcommand{\maps}{\colon}

\newcommand{\id}{{\rm id}}

\newcommand{\Res}{{\rm Res}}
\newcommand{\End}{{\rm End}}
\newcommand{\chr}{{\rm ch}}

\newcommand{\scs}{\scriptstyle}

\theoremstyle{definition}
\newtheorem{thm}{Theorem}[section]
\newtheorem{cor}[thm]{Corollary}
\newtheorem{conj}[thm]{Conjecture}
\newtheorem{lem}[thm]{Lemma}
\newtheorem{rem}[thm]{Remark}
\newtheorem{prop}[thm]{Proposition}

\newtheorem{example}[thm]{Example}

        \newcommand{\be}{\begin{equation}}
        \newcommand{\ee}{\end{equation}}
        \newcommand{\ba}{\begin{eqnarray}}
        \newcommand{\ea}{\end{eqnarray}}
        \newcommand{\ban}{\begin{eqnarray*}}
        \newcommand{\ean}{\end{eqnarray*}}
        \newcommand{\barr}{\begin{array}}
        \newcommand{\earr}{\end{array}}

\numberwithin{equation}{section}

\def\emph#1{{\sl #1\/}}

\let\hat=\widehat
\let\tilde=\widetilde

\let\phi=\varphi
\let\epsilon=\varepsilon

\usepackage{bbm}
\def\C{{\mathbbm C}}
\def\N{{\mathbbm N}}
\def\R{{\mathbbm R}}
\def\Z{{\mathbbm Z}}
\def\Q{{\mathbbm Q}}

\def\cal#1{\mathcal{#1}}%
\def\1{\mathbbm{1}}%
\def\nn{\notag}

\def\Res{{\mathrm{Res}}}
\def\Ind{{\mathrm{Ind}}}
\def\lra{{\longrightarrow}}
\def\dmod{{\mathrm{-mod}}}   
\def\fmod{{\mathrm{-fmod}}}   
\def\pmod{{\mathrm{-pmod}}}  
\def\gdim{{\mathrm{gdim}}}
\def\Ext{{\mathrm{Ext}}}
\def\pseq{{\mathrm{Seqd}}}
\def\Id{\mathrm{Id}}
\def\mc{\mathcal}
\def\mf{\mathfrak}
\def\Af{{_{\mc{A}}\mathbf{f}}}    
\def\shuffle{\,\raise 1pt\hbox{$\scriptscriptstyle\cup{\mskip
               -4mu}\cup$}\,}
\newcommand{\define}{\stackrel{\mbox{\scriptsize{def}}}{=}}

\newcommand{\up}[1]{\xybox{
   (-3,-13)*{};
  (3,8)*{};
 (0,0)*{\includegraphics[scale=0.5]{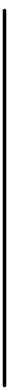}};
 (-.1,-12)*{\scs #1};
 }}

\renewcommand{\sup}[1]{\xybox{
   (-3,-7)*{};
  (3,6)*{};
 (0,0)*{\includegraphics[scale=0.5]{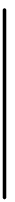}};
 (-.1,-7)*{\scs #1};
 }}

 \newcommand{\supdot}[1]{\xybox{
   (-3,-7)*{};
  (3,6)*{};
 (0,0)*{\includegraphics[scale=0.5]{short_up.eps}};
 (-.1,-7)*{\scs #1}; (0,0)*{\bullet};
 }}

\newcommand{\dcross}[2]{\xybox{
 (-6,-7)*{};
 (6,6)*{};
 (0,0)*{\includegraphics[scale=0.5]{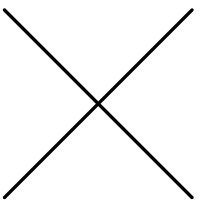}};
 (-5.1,-7)*{\scs #1};
 (4.7,-7)*{\scs #2};
}}

\newcommand{\ddcross}[2]{\xybox{
(-2.5,-2.5)*{\bullet}; (2.5,-2.5)*{\bullet};
 (-6,-7)*{};
 (6,6)*{};
 (0,0)*{\includegraphics[scale=0.5]{dcross.eps}};
 (-5.1,-7)*{\scs #1};
 (4.7,-7)*{\scs #2};
}}

\newcommand{\dcrossul}[2]{\xybox{
(-2.5,2.5)*{\bullet};
 (-6,-7)*{};
 (6,6)*{};
 (0,0)*{\includegraphics[scale=0.5]{dcross.eps}};
 (-5.1,-7)*{\scs #1};
 (4.7,-7)*{\scs #2};
}}

\newcommand{\dcrossur}[2]{\xybox{
(2.5,2.5)*{\bullet};
 (-6,-7)*{};
 (6,6)*{};
 (0,0)*{\includegraphics[scale=0.5]{dcross.eps}};
 (-5.1,-7)*{\scs #1};
 (4.7,-7)*{\scs #2};
}}

\newcommand{\dcrossdl}[2]{\xybox{
(-2.5,-2.5)*{\bullet};
 (-6,-7)*{};
 (6,6)*{};
 (0,0)*{\includegraphics[scale=0.5]{dcross.eps}};
 (-5.1,-7)*{\scs #1};
 (4.7,-7)*{\scs #2};
}}

\newcommand{\dcrossdr}[2]{\xybox{
(2.5,-2.5)*{\bullet};
 (-6,-7)*{};
 (6,6)*{};
 (0,0)*{\includegraphics[scale=0.5]{dcross.eps}};
 (-5.1,-7)*{\scs #1};
 (4.7,-7)*{\scs #2};
}}

\newcommand{\twocross}[2]{\xybox{
 (-6,-13)*{};
 (6,8)*{};
 (0,0)*{\includegraphics[scale=0.5]{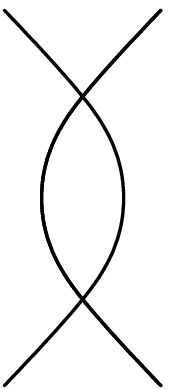}};
 (-4.1,-12)*{\scs #1};
 (3.7,-12)*{\scs #2};
}}

\newcommand{\linecrossL}[3]{\xybox{
 (-6,-13)*{};
 (6,8)*{};
 (0,0)*{\includegraphics[scale=0.5]{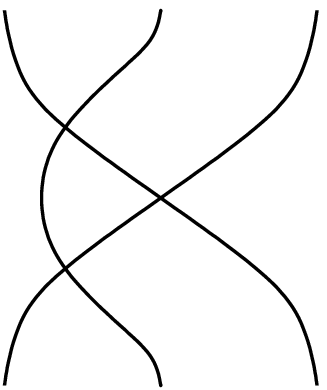}};
 (-8,-12)*{\scs #1};
 (0,-12)*{\scs #2};
 (8,-12)*{\scs #3};
}}

\newcommand{\linecrossR}[3]{\xybox{
 (-6,-13)*{};
 (6,8)*{};
 (0,0)*{\includegraphics[angle=180, scale=0.5]{line_crossL.eps}};
 (-8,-12)*{\scs #1};
 (0,-12)*{\scs #2};
 (8,-12)*{\scs #3};
}}

\newcommand{\eqdefB}{
B_0=\left(
 \begin{array}{c}
   \xy (-6,-13)*{};
 (6,8)*{};
 (-5.2,5)*{\scs \bullet};
 (0,0)*{\reflectbox{\includegraphics[scale=0.35]{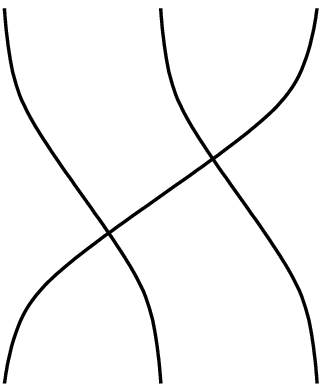}}};
 (-6,-9)*{\scs i};(0,-9)*{\scs j};(6,-9)*{\scs i};
 \endxy
  \\
  \xy (-6,-13)*{};
 (6,8)*{};
 (0.1,5)*{\scs \bullet};
 (0,0)*{\includegraphics[scale=0.35]{j1i.eps}};
 (-6,-9)*{\scs i};(0,-9)*{\scs j};(6,-9)*{\scs i};
 \endxy
 \end{array}
\right), \qquad \quad B_1=\left(
 \begin{array}{ccc}
 -\xy (-6,-13)*{};
 (6,8)*{};
 (0,0)*{\includegraphics[scale=0.35]{j1i.eps}};
 (-6,-9)*{\scs i};(0,-9)*{\scs i};(6,-9)*{\scs j};
 \endxy
    & &
 \xy (-6,-13)*{};
 (6,8)*{};
 (0,0)*{\reflectbox{\includegraphics[scale=0.35]{j1i.eps}}};
 (-6,-9)*{\scs j};(0,-9)*{\scs i};(6,-9)*{\scs i};
 \endxy
 \end{array}
\right) }

\newcommand{\BoBoneA}{
  B_0 B_1 =  \left(
  \begin{array}{c}
   \xy
     (-6,-9)*{}; (6,9)*{};
    (-3,2)*{\scs \bullet};
    (0,0)*{\reflectbox{\includegraphics[scale=0.25]{j1i.eps}}};
    (-4,-7)*{\scriptscriptstyle i};
    (0,-7)*{\scriptscriptstyle j};
    (4,-7)*{\scriptscriptstyle i};
    \endxy
     \\
    \vcenter{\xy (-6,-9)*{}; (6,9)*{};
    (0.3,2)*{\scs \bullet};
    (0,0)*{\includegraphics[scale=0.25]{j1i.eps}};
    (-4,-7)*{\scriptscriptstyle i};
    (0,-7)*{\scriptscriptstyle j};
    (4,-7)*{\scriptscriptstyle i};
    \endxy}
 \end{array}
 \right)
 \left(
 \begin{array}{ccc}
 -\vcenter{\xy  (-6,-9)*{}; (6,8)*{};
 (0,0)*{\includegraphics[scale=0.25]{j1i.eps}};
  (-4,-7)*{\scriptscriptstyle i};
 (0,-7)*{\scriptscriptstyle i};
 (4,-7)*{\scriptscriptstyle j};
 \endxy}
    & &
 \vcenter{\xy (-6,-9)*{}; (6,8)*{};
 (0,0)*{\reflectbox{\includegraphics[scale=0.25]{j1i.eps}}};
  (-4,-7)*{\scriptscriptstyle j};
 (0,-7)*{\scriptscriptstyle i};
 (4,-7)*{\scriptscriptstyle i};
 \endxy}
 \end{array}
\right) \;  = \; \left(
\begin{array}{cc}
 - \vcenter{\xy (-6,-9)*{}; (6,15)*{};
 (-3.1,12)*{\scs \bullet};
 (0,0)*{\includegraphics[scale=0.25]{j1i.eps}};
 (0,9.5)*{\reflectbox{\includegraphics[angle=180, scale=0.25]{j1i.eps}}};
 (-4,-7)*{\scriptscriptstyle i};
 (0,-7)*{\scriptscriptstyle i};
 (4,-7)*{\scriptscriptstyle j};
 \endxy}  &
 \vcenter{\xy (-6,-9)*{}; (6,15)*{};
 (-3.1,12)*{\scs \bullet};
 (0,0)*{\reflectbox{\includegraphics[scale=0.25]{j1i.eps}}};
 (0,9.5)*{\reflectbox{\includegraphics[angle=180, scale=0.25]{j1i.eps}}};
 (-4,-7)*{\scriptscriptstyle j};
 (0,-7)*{\scriptscriptstyle i};
 (4,-7)*{\scriptscriptstyle i};
 \endxy}\\
  - \vcenter{\xy (-6,-9)*{}; (6,16)*{};
 (.3,12)*{\scs \bullet};
 (0,0)*{\includegraphics[scale=0.25]{j1i.eps}};
 (0,9.5)*{\includegraphics[angle=180, scale=0.25]{j1i.eps}};
 (-4,-7)*{\scriptscriptstyle i};
 (0,-7)*{\scriptscriptstyle i};
 (4,-7)*{\scriptscriptstyle j};
 \endxy }  &
 \vcenter{ \xy (-6,-9)*{}; (6,16)*{};
 (.3,12)*{\scs \bullet};
 (0,0)*{\reflectbox{\includegraphics[scale=0.25]{j1i.eps}}};
 (0,9.5)*{\includegraphics[angle=180, scale=0.25]{j1i.eps}};
 (-4,-7)*{\scriptscriptstyle j};
 (0,-7)*{\scriptscriptstyle i};
 (4,-7)*{\scriptscriptstyle i};
 \endxy}
\end{array}
\right) }

\newcommand{\BoBoneB}{
\;\; \refequal{\eqref{eq_UUzero}} \;\;
 \left(
 \begin{array}{cc}
    - \xy (-6,-9)*{}; (6,8)*{};
    (-1.2,4.5)*{\scs \bullet}; (1.2,0)*{\scs \bullet};
     (0,0)*{\includegraphics[scale=0.3]{R2.eps}};
     (6,.5)*{\includegraphics[scale=0.33]{up.eps}};
      (-2.5,-8)*{\scriptscriptstyle i};
    (2.5,-8)*{\scriptscriptstyle i};
    (6,-8)*{\scriptscriptstyle j};
    \endxy
    & 0 \\
   0 & \xy (-6,-9)*{}; (6,8)*{};
   (2.8,4.5)*{\scs \bullet}; (2.8,0)*{\scs \bullet};
     (4,0)*{\includegraphics[scale=0.3]{R2.eps}};
     (-2,.5)*{\includegraphics[scale=0.33]{up.eps}};
      (-2.5,-8)*{\scriptscriptstyle j};
    (2.5,-8)*{\scriptscriptstyle i};
    (6,-8)*{\scriptscriptstyle i};
    \endxy
 \end{array}
\right) \nn \;\; \refequal{\eqref{eq_UUzero},\eqref{eq_iislide1}}\;\; \left(
 \begin{array}{cc}
  \vcenter{ \xy
   (-5,-7)*{}; (-1.5,1.5)*{\scs \bullet};
   (7,0)*{\includegraphics[scale=0.35]{short_up.eps}};
   (0,0.1)*{\includegraphics[scale=0.3]{dcross.eps}};
   (-3,-5)*{\scriptscriptstyle i};
    (3,-5)*{\scriptscriptstyle i};
    (7,-5)*{\scriptscriptstyle j};
   \endxy }& 0 \\
   0 &
   \vcenter{\xy
   (-5,-7)*{}; (2.2,1.5)*{\scs \bullet};
   (-3,0)*{\includegraphics[scale=0.35]{short_up.eps}};
   (4,0.1)*{\includegraphics[scale=0.3]{dcross.eps}};
   (-3,-5)*{\scriptscriptstyle j};
    (1,-5)*{\scriptscriptstyle i};
    (7,-5)*{\scriptscriptstyle i};
   \endxy}
 \end{array}
\right) }

\newcommand{\BoneBoA}{
B_1 B_0 = \left(
 \begin{array}{ccc}
 -\vcenter{\xy  (-6,-9)*{}; (6,8)*{};
 (0,0)*{\includegraphics[scale=0.25]{j1i.eps}};
  (-4,-7)*{\scriptscriptstyle i};
 (0,-7)*{\scriptscriptstyle i};
 (4,-7)*{\scriptscriptstyle j};
 \endxy}
    & &
 \vcenter{\xy (-6,-9)*{}; (6,8)*{};
 (0,0)*{\reflectbox{\includegraphics[scale=0.25]{j1i.eps}}};
  (-4,-7)*{\scriptscriptstyle j};
 (0,-7)*{\scriptscriptstyle i};
 (4,-7)*{\scriptscriptstyle i};
 \endxy}
 \end{array}
\right)  \left(
 \begin{array}{c}
   \xy
 (-6,-9)*{}; (6,9)*{};
 (-3,2)*{\scs \bullet};
 (0,0)*{\reflectbox{\includegraphics[scale=0.25]{j1i.eps}}};
 (-4,-7)*{\scriptscriptstyle i};
 (0,-7)*{\scriptscriptstyle j};
 (4,-7)*{\scriptscriptstyle i};
 \endxy
  \\
  \vcenter{\xy (-6,-9)*{}; (6,9)*{};
 (0.3,2)*{\scs \bullet};
 (0,0)*{\includegraphics[scale=0.25]{j1i.eps}};
 (-4,-7)*{\scriptscriptstyle i};
 (0,-7)*{\scriptscriptstyle j};
 (4,-7)*{\scriptscriptstyle i};
 \endxy}
 \end{array}
\right) \;\;=\;\; - \vcenter{ \xy (-6,-9)*{}; (6,16)*{};
 (-4,4)*{\scs \bullet};
 (0,0)*{\reflectbox{\includegraphics[scale=0.25]{j1i.eps}}};
 (0,9.5)*{\includegraphics[angle=180, scale=0.25]{j1i.eps}};
 (-4,-7)*{\scriptscriptstyle i};
 (0,-7)*{\scriptscriptstyle j};
 (4,-7)*{\scriptscriptstyle i};
 \endxy}
 \;\;+\;\;
\vcenter{\xy (-6,-9)*{}; (6,15)*{};
 (4,4)*{\scs \bullet};
 (0,0)*{\includegraphics[scale=0.25]{j1i.eps}};
 (0,9.5)*{\reflectbox{\includegraphics[angle=180, scale=0.25]{j1i.eps}}};
 (-4,-7)*{\scriptscriptstyle i};
 (0,-7)*{\scriptscriptstyle j};
 (4,-7)*{\scriptscriptstyle i};
 \endxy}
 }

\newcommand{\BoneBoB}{
  \;\;\refequal{\eqref{eq_UUzero}}\;\;
  -\xy
  (0,0)*{\includegraphics[angle=180, scale=0.3]{line_crossL.eps}};
   (-5,-8)*{\scriptscriptstyle i};
 (0,-8)*{\scriptscriptstyle j};
 (5,-8)*{\scriptscriptstyle i};
  \endxy
\;\; + \;\;
  \xy
  (0,0)*{\includegraphics[scale=0.3]{line_crossL.eps}};
   (-5,-8)*{\scriptscriptstyle i};
 (0,-8)*{\scriptscriptstyle j};
 (5,-8)*{\scriptscriptstyle i};
  \endxy
\;\;\refequal{\eqref{eq_r3_hard}}\;\;
 \xy
   (-5,0)*{\includegraphics[scale=0.3]{up.eps}};
   (0,0)*{\includegraphics[scale=0.3]{up.eps}};
   (5,0)*{\includegraphics[scale=0.3]{up.eps}};
   (-5,-8)*{\scriptscriptstyle i};
 (0,-8)*{\scriptscriptstyle j};
 (5,-8)*{\scriptscriptstyle i};
  \endxy
}

\title{A diagrammatic approach to categorification of quantum groups I}
      \author{ Mikhail Khovanov and Aaron D.\ Lauda}

\begin{document}

\date{March 31, 2008}

\maketitle

\begin{abstract}
To each graph without loops and multiple edges we assign a family of rings.
Categories of projective modules over these rings categorify $U^-_q(\mathfrak{g})$,
where $\mathfrak{g}$ is the Kac-Moody Lie algebra associated with the graph.
\end{abstract}

\tableofcontents

%
\section{Introduction}
%

The goal of this paper is to categorify $U^-=U_q^-(\mf{g})$, for an arbitrary
simply-laced Kac-Moody algebra $\mf{g}$. Here $U^-$ stands for the quantum
deformation of the universal enveloping algebra of the ``lower-triangular'' subalgebra
of $\mf{g}$.

Following the discovery of quantum groups $U_q(\mf{g})$ by Drinfeld~\cite{Drinfeld} and
Jimbo~\cite{Jimbo}, Ringel~\cite{Ringel} found a Hall algebra interpretation of
the negative half $U^-$ of the quantum group in the simply-laced Dynkin case.
Lusztig~\cite{Lus1}, \cite{Lus2}, \cite{Lus3} gave a geometric interpretation
of $U^-$ and produced a canonical basis there via a sophisticated approach which
required the full strength of the theory of $l$-adic perverse sheaves. Kashiwara~\cite{Kas1}
defined a crystal basis of $U^-$ at $0$, a graph equipped with extra data, and in~\cite{Kas2}
constructed the so-called global crystal basis of $U^-$. Grojnowski and Lusztig~\cite{GL2}
proved that the global crystal basis and the canonical basis are the same. The canonical basis
$\mathbf{B}$ of $U^-$ gives rise to bases in all irreducible integrable $U$-representations.
Lusztig~\cite{Lus4} also produced an idempotented version $\U$ of $U$ and defined a basis
there.

The work of Ariki~\cite{Ari1} can be viewed as a categorification of the restricted
dual of $U^-(\mf{g})$ for $\mf{g}=\mf{sl}_N$ and $\mf{g}=\widehat{\mf{sl}}_N$
and a categorification of all irreducible integrable representations of these Lie algebras
(see also~\cite{LLT}, \cite{Ari2}, \cite{Ari3}, \cite{Mathas}).
An integral version of the restricted dual of  $U^-(\mf{g})$ becomes the sum of
Grothendieck groups of suitable blocks of affine Hecke algebra representations.
An earlier work of Zelevinsky~\cite{Zel}
can be understood in this context as a parametrization of basis elements of
$U^-(\mf{g})^{\ast}$
via certain irreducible representations of affine Hecke algebras. Irreducible integrable
representations of $U(\mf{g})$ become Grothendieck groups of Ariki-Koike cyclotomic Hecke
algebras, which are certain finite-dimensional quotient algebras of affine Hecke algebras.

Grojnowski~\cite{Groj} found a purely algebraic way to understand these categorifications
via a generalization of Kleshchev's methods for studying modular representations of
the symmetric group~\cite{Kle1}, \cite{Kle2}, \cite{Kle3}. This approach was
further developed by Grojnowski-Vazirani~\cite{GV}, Vazirani~\cite{Vaz1}, \cite{Vaz2},
Brundan-Kleshchev~\cite{BK} and others. It is explained in Kleshchev~\cite{KleBook} in
the context of degenerate affine Hecke algebras.

In this paper we introduce graded algebras categorifying $U^-_q(\mf{g})$, for
an arbitrary simply-laced $\mf{g}$. We start with an unoriented graph
$\Gamma$ without loops and multiple edges. Let $I$ be the set of vertices of $\Gamma$.
The bilinear Cartan form on $\N[I]$ is given on the basis elements $i,j\in I$ by
$$
      i\cdot j = \begin{cases} 2 & \textrm{if $i=j$}, \\
      -1 & \textrm{if $i$ and $j$ are joined by an edge}, \\
      0 & \textrm{otherwise}. \end{cases}
$$
The algebra
$U^-$ over $\Q(q)$, the negative (or positive) half of the quantum universal enveloping
algebra, has generators $\theta_i$, $i\in I$, and defining relations
\begin{eqnarray*}
  \theta_i \theta_j & = & \theta_j \theta_i \  \   \mathrm{if} \  \   i\cdot j = 0, \\
  (q+q^{-1})\theta_i \theta_j \theta_i & = & \theta_i^2 \theta_j + \theta_j \theta_i^2
    \  \   \mathrm{if} \ \   i\cdot j = -1.
\end{eqnarray*}
The algebra $U^-$ contains a subring $\Af$, which is the $\Z[q,q^{-1}]$-lattice generated by
all products of quantum divided powers $\theta_i^{(a)}$. The canonical basis $\mathbf{B}$
is a basis of $\Af$ viewed as a free  $\Z[q,q^{-1}]$-module.

In Section~\ref{sec-ringsrnu} of this paper to each graph $\Gamma$ as above
we assign a family of graded rings $R(\nu)$, over $\nu\in \N[I]$.
The rings are defined geometrically, via
braid-like plane diagrams which consist of interacting strings labelled by vertices of the graph.
We prove basic results about these rings, then switch from the ground ring $\Z$ to a field
$\Bbbk$ to simplify the study of $R(\nu)$-modules. We show that
the representation theory of $R(\nu)$ categorifies the
integral form $\Af$ of $U^-$. We consider
the category $R(\nu)\pmod$ of finitely-generated graded left projective $R(\nu)$-modules
and its Grothendieck group $K_0(R(\nu))$. Let $R=\oplus_{\nu} R(\nu)$ and
define
$$K_0(R)= \bigoplus_{\nu\in \N[I]} K_0(R(\nu)).$$
Induction and restriction functors coming from
the  inclusions $R(\nu)\otimes R(\nu')\subset R(\nu+\nu')$ give rise to the multiplication
and comultiplication homomorphisms
 $$ K_0(R) \otimes K_0(R) \longrightarrow K_0(R), \ \ K_0(R)\longrightarrow
   K_0(R) \otimes K_0(R) $$
that satisfy the same properties as those for $\Af$. We define a homomorphism
of $\Z[q,q^{-1}]$-algebras
$\gamma: \Af \lra K_0(R)$ that also respects comultiplication and takes a divided powers
product element $\theta=\theta_{i_1}^{(a_1)}\dots \theta_{i_r}^{(a_r)}$ to the
image of a certain projective module $P_{\theta}$ in the Grothendieck group.

The quantum Gabber-Kac theorem implies that $\gamma$ is injective.
By mirroring for the case of rings $R(\nu)$ the methods of
Kleshchev, Grojnowski, and Vazirani, who studied socles of
induction and restriction applied to irreducible representations,
we show that homomorphism $\gamma$ is surjective for any graph $\Gamma$ and
any field $\Bbbk$. The main result of the paper is the following theorem.

\begin{thm} \label{thm_first}
 $\gamma: \Af\lra K_0(R)$ is an isomorphism of $\N[I]$-graded twisted
bialgebras.
\end{thm}

The term ``twisted bialgebras'' is used above, since the comultiplication in $\Af$ and
$K_0(R)$ becomes an algebra homomorphism only after the multiplication in the
tensor squares $\Af^{\otimes 2}$ and $K_0(R)^{\otimes 2}$ is twisted by powers
of $q$.

We conjecture that, when $\Gamma$ is a tree and $\Bbbk=\C$, isomorphism $\gamma$
takes canonical basis elements to the images of indecomposable projective modules in $K_0(R)$.
When the graph is a single vertex, this
conjecture  is almost trivial. We verify the conjecture in a simple case of the graph
$\Gamma$ with two vertices and one edge, with all canonical basis elements being monomials.

The rings $R(\nu)$ should be linked to Lusztig's geometric realization of $\Af$:

\begin{conj}\label{conj-perv}
For $\Bbbk=\C$ and graph $\Gamma$ a tree, the
algebra $R(\nu)$ is Morita equivalent to the algebra of equivariant ext groups
$\Ext_{G_{\nu}}(L,L)$, where $G_{\nu}=\prod_{i\in I} GL(\nu_i)$ and
$L$ is the sum of simple perverse sheaves $L_b$, over all $b\in \mathbf{B}_{\nu}$,
in Lusztig's geometric realization of $\Af$.
\end{conj}
This conjecture should follow from an isomorphism between $R(\nu)$ and a
suitable convolution algebra. When $\Gamma$ contains cycles, it's possible to
modify $R(\nu)$ by introducing ``monodromies'' around the cycles, and the conjecture
is likely to hold for a modified version of $R(\nu)$.

Our results hint at the relation between representation theory of affine Hecke
algebras for $GL(n)$ when $q$ is not a root of unity and
representations of $R(\nu)$ when the graph $\Gamma$ is a chain. We conjecture that
completions of affine Hecke algebras along suitable maximal central ideals are
Morita equivalent (or even isomorphic) to completions of  $R(\nu)$ along the grading.
The above conjectures, if true, would link Lusztig's geometrization of $U^-$ with
Ariki's categorification of $U^-$ and its restricted dual for $q=1$ and $\Gamma=A_n$.

We arrived at the definition of rings $R(\nu)$ from computations involving
homomorphisms of bimodules over cohomology rings of partial flag varieties.
The bimodules themselves are the cohomology groups of partial and iterated
flag varieties that give correspondences for the action of generators $e_i$ and $f_i$ of
quantum $\mathfrak{sl}_N$ in the Beilinson-Lusztig-MacPherson geometric
model~\cite{BLM} of quantum $\mathfrak{sl}_N$. This model was given
a categorical interpretation by Grojnowski and Lusztig~\cite{GL1} and later
reinterpreted, for $N=2$, via translation and Zuckerman functors in~\cite{BFK},
with various generalizations constructed in~\cite{FKS},~\cite{Sus},~\cite{Zheng}, 
and in a very recent striking work~\cite{Zheng2}. 

In Section~\ref{subsec-conj-irr} we define certain quotient algebras of $R(\nu)$
and conjecture that their categories of modules categorify irreducible
integrable representations of $U_q(\mf{g})$. A straightforward generalization
of our constructions and results
from algebras $R(\nu)$ and their quotient algebras in the simply-laced case to that of an
arbitrary symmetrizable Kac-Moody algebra $\mf{g}$ will be presented in a follow-up paper.

In the paper~\cite{CR} (see also~\cite{Rou1}), Chuang and Rouquier defined
$sl(2)$-categorifications, substantiated them with many diverse examples,
and applied to the modular representation theory.
Partially inspired by~\cite{CR} and \cite{CF}, the second author
suggested and investigated a categorification of Lusztig's idempotent completion
$\U$ of quantum $sl(2)$ in~\cite{Lau1}, \cite{Lau2}.
A definition of a categorification of $\U_q(\mf{g})$ for any simply-laced $\mf{g}$ can
be obtained by combining the diagrammatic relations of $R(\nu)$ with those of
$\U$-categorification in~\cite{Lau1}.

R.~Rouquier~\cite{Rou2}, in his recent talk,
defined $sl(N)$- and affine $sl(N)$-categorifications and outlined a conjectural program 
that aims
to vastly generalize his prior work with J.~Chuang~\cite{CR} on $sl(2)$-categorifications.
We expect Rouquier's and our approaches to be closely related. R.~Rouquier informed 
us that a signed version of rings $R(\nu)$ appears in his categorification~\cite{Rou3} 
of $U(\mf{g})$ for a simply-laced $\mf{g}$. 

\vspace{0.1in}

{\bf Acknowledgments:} This paper was written while the first author was
at the Institute for Advanced Study, and he would like to thank the Institute for its
hospitality and the NSF for fully supporting him during that time through the
grant DMS-0635607. Additional partial support came from the NSF grant DMS-0706924.

%
\section{Rings $R(\nu)$ and their properties}
\label{sec-ringsrnu}
%

%
\subsection{Definitions}
\label{subsec-definitions}
%

We fix a graph $\Gamma$, not necessarily finite,  with set of vertices $I$ and
unoriented edges $E_{\Gamma}$. We require that $\Gamma$ has no loops
and multiple edges. By $\N[I]$ we denote the commutative semigroup freely
generated by vertices of $\Gamma$ and for $\nu \in \N[I]$ write
\begin{equation}
  \nu = \sum_{i \in I} \nu_{i} \cdot i \;, \ \  \nu_i \in \N, \ \
  \N=\{0,1,2,\dots\}.
\end{equation}
Let $|\nu|=\sum \nu_i \in \N$, and ${\rm Supp}(\nu)=\{ i \mid \nu_i \neq 0\}$.
We define a bilinear form on $\Z[I]$ by $i \cdot i=2$, $i \cdot j=-1$ if $i$ and
$j$ are connected by an edge, and $i \cdot j=0$ otherwise. In the basis $\{i\}_{i
\in I}$ of vertices the bilinear form is given by the Cartan matrix of $\Gamma$.

To $\Gamma$ we associate a diagrammatic calculus of planar diagrams. We
consider collections of arcs on the plane connecting $m$ points on one horizontal
line with $m$ points on another horizontal line.  The position of $m$ points on
the line is fixed once and for all (for instance, we could take points $\{ 1,2,
\dots, m \} \in \R$).  Arcs have no critical points when projected to the
$y$-axis of the plane (a condition reminiscent of braids).  Each arc is labelled
by a vertex of $\Gamma$. Arcs can intersect, but no triple intersections are
allowed. An arc can carry dots.  An example of such a diagram is shown
below,
\begin{equation} \label{eq_example_diagram} \xy
  (-12,-13)*{};
  (0,0)*{\includegraphics[scale=0.5]{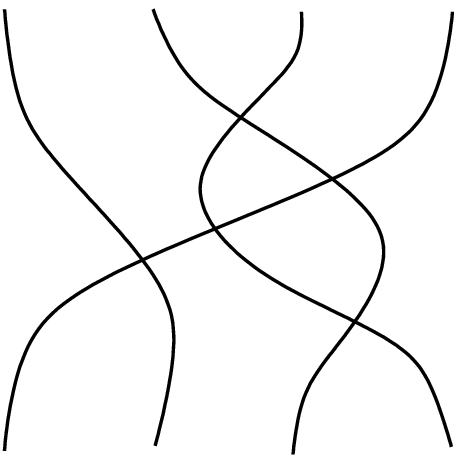}};
  (-12,-14)*{\scs i};  (-4,-14)*{\scs j};
  (4,-14)*{\scs i};   (11.2,-14)*{\scs k};
  (-10.7,-8)*{\bullet};(10.7,7)*{\bullet};
  (-1.2,3)*{\bullet};(8,-1)*{\bullet};
 \endxy \end{equation}
where $i,j,k$ are vertices of $\Gamma$ and the label of an arc is written at the
bottom end of the arc. We allow isotopies that do not change the combinatorial
type of the diagram and do not create critical points for the projection onto the
$z$-axis.

\[
 \xy (0,0)*{\includegraphics[scale=0.5]{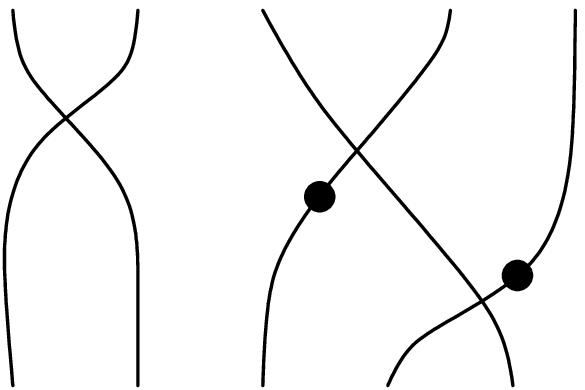}};
 (-15,-12)*{\scs i};  (-8,-12)*{\scs j};
  (-2,-12)*{\scs i};   (12,-12)*{\scs k};(5,-12)*{\scs k};\endxy
 \qquad \sim \qquad
 \xy (0,0)*{\includegraphics[scale=0.5]{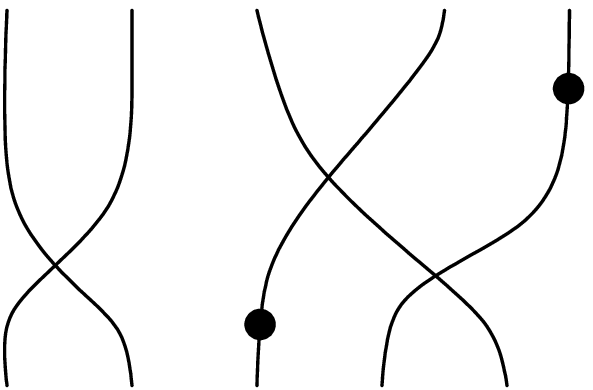}};
 (-15,-12)*{\scs i};  (-8,-12)*{\scs j};
  (-2,-12)*{\scs i};   (12,-12)*{\scs k};(5,-12)*{\scs k};  \endxy
\]
We proceed by allowing finite linear combinations of these diagrams with integral
coefficients, modulo the following local relations
\begin{eqnarray} \label{eq_UUzero}
   \xy   (0,0)*{\twocross{i}{j}}; \endxy
 & = & \left\{
\begin{array}{ccc}
  0 & \qquad & \text{if $i=j$, } \\ \\
  \xy (0,0)*{\sup{i}};  (8,0)*{\sup{j}};  \endxy
  & &
 \text{if $i \cdot j=0$, }
  \\    \\
  \xy  (0,0)*{\supdot{i}};   (8,0)*{\sup{j}};  \endxy
  \quad + \quad
   \xy  (8,0)*{\supdot{j}};  (0,0)*{\sup{i}};  \endxy
 & &
 \text{if $i \cdot j=-1$. }
\end{array}
\right.
\end{eqnarray}

\begin{eqnarray}\label{eq_ijslide}
  \xy  (0,0)*{\dcrossul{i}{j}};  \endxy
 \quad  = \;\;
   \xy  (0,0)*{\dcrossdr{i}{j}};   \endxy
& \quad &
   \xy  (0,0)*{\dcrossur{i}{j}};  \endxy
 \quad = \;\;
   \xy  (0,0)*{\dcrossdl{i}{j}};  \endxy \qquad \text{for $i \neq j$}
\end{eqnarray}

\begin{eqnarray}        \label{eq_iislide1}
 \xy  (0,0)*{\dcrossul{i}{i}}; \endxy
    \quad - \quad
 \xy (0,0)*{\dcrossdr{i}{i}}; \endxy
  & = &
 \xy (-3,0)*{\sup{i}}; (3,0)*{\sup{i}}; \endxy \\      \label{eq_iislide2}
  \xy (0,0)*{\dcrossdl{i}{i}}; \endxy
 \quad - \quad
 \xy (0,0)*{\dcrossur{i}{i}}; \endxy
  & = &
 \xy (-3,0)*{\sup{i}}; (3,0)*{\sup{i}}; \endxy
\end{eqnarray}

\begin{eqnarray}      \label{eq_r3_easy}
\xy  (0,0)*{\linecrossL{i}{j}{k}}; \endxy
  &=&
\xy (0,0)*{\linecrossR{i}{j}{k}}; \endxy
 \qquad \text{unless $i=k$ and $i \cdot j=-1$   \hspace{1in} }
\\                   \label{eq_r3_hard}
\xy (0,0)*{\linecrossL{i}{j}{i}}; \endxy
  &-&
\xy (0,0)*{\linecrossR{i}{j}{i}}; \endxy
 \quad = \quad
 \xy  (-9,0)*{\up{i}};  (0,0)*{\up{j}}; (9,0)*{\up{i}}; \endxy
 \qquad \text{if $i \cdot j=-1$ }
\end{eqnarray}

Fix $\nu \in \N[I]$.  Let $\seq(\nu)$ be the set of all sequences of vertices $\ii = i_1
\ldots i_m$ where $i_k \in I$ for each $k$ and vertex $i$ appears $\nu_i$ times
in the sequence.  The length $m$ of the sequence is equal to $|\nu|$ and the
cardinality of $\seq(\nu)$ is equal to $ \sbinom{\nu}{\nu_i,\nu_j, \ldots}$,
taken over all $i \in I$.  For instance,
$$\seq(2i+j)=\{ iij,iji,jii\}.$$  Each diagram $D$ as described above determines
two sequences $\bot(D)$ and $\top(D)$ in $\seq(\nu)$ for some $\nu$. The sequence
$\bot(D)$ is given by reading the labels of arcs of $D$ at the bottom position
from left to right.  We define $\top(D)$ likewise. For instance, for the diagram
in \eqref{eq_example_diagram}, $\bot(D)=ijik$ and $\top(D)=jiki$. We often
abbreviate sequences with many equal consecutive terms, and write
$i_1^{n_1}\dots i_r^{n_r}$ for $i_1\dots i_1 i_2\dots i_2 \dots i_r\dots i_r$,
where $n_1+\dots + n_r =m$.

Define the ring $R(\nu)$ as follows:
\begin{equation}
  R(\nu) = \bigoplus_{\ii,\jj \in \seq(\nu)} \jRi
\end{equation}
as an abelian group, where $\jRi$ is the abelian group of all linear combinations
of diagrams with $\bot(D)=\ii$ and $\top(D)=\jj$ modulo the relations
\eqref{eq_UUzero}--\eqref{eq_r3_hard}. The product in $R(\nu)$ is given by
concatenation (see the left diagram in \eqref{eq_two_diag} below)
\begin{equation}
 {_{\kk}R(\nu)_{\jj}} \otimes \jRi \to {_{\kk}R(\nu)_{\ii}},
\end{equation}
and $xy=0$ for $x \in {_{\textbf{\textit{l}}}R(\nu)_{\kk}}$ and $y \in \jRi$ if
$\kk \neq \jj$.

\begin{equation} \label{eq_two_diag}
  \xy  (0,0)*{\includegraphics[scale=0.6]{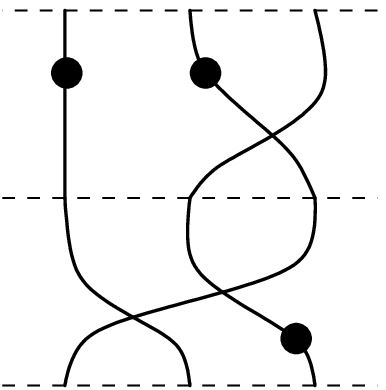}};
  (14,11)*{\scs \kk}; (14,0)*{\scs \jj}; (14,-11)*{\scs \ii}; \endxy
 \qquad \qquad
 \xy (0,18)*{}; (0,0)*{\up{\;\;i_1}}; (4,0)*{\up{\;\;i_2}};
  (8,0)*{\up{ }}; (12,0)*{\dots}; (16,0)*{\up{\;\;\;i_m}};
 \endxy\end{equation}

By construction, $R(\nu)$ is an associative ring.  For each $\ii
\in \seq(\nu)$ the diagram $1_{\ii} \in \iRi$ shown on the right of
\eqref{eq_two_diag} is an idempotent, $1_{\ii}^2=1_{\ii}$, $x1_{\ii}=x$ for all
$x \in \jRi$ and $1_{\ii}x=x$ for all $x \in \iRj$, for all $\jj$. Furthermore,
$1 = \sum_{\ii \in \seq(\nu)}1_{\ii}$ is the unit element of $R(\nu)$.
We turn $R(\nu)$ into a graded ring by declaring degrees of the generators to be
\begin{equation}
  \deg\left( \xy  (-3,0)*{\supdot{i}}; \endxy \right) = 2, \qquad \deg\left(  \xy
  (0,0)*{\dcross{i}{j}};  \endxy \right) = - i \cdot j.
\end{equation}
Let
\begin{equation}
  P_{\ii}  =  \bigoplus_{\jj \in \seq(\nu)} \jRi ,  \  \
  {_{\jj}P} = \bigoplus_{\ii \in \seq(\nu)} \jRi.
\end{equation}
$P_{\ii}$ is a left graded projective $R(\nu)$-module and ${_{\jj}P}$ is a right
graded projective $R(\nu)$-module.

Flipping a diagram about a horizontal axis induces a grading--preserving
antiinvolution $\psi$ of $R(\nu)$ which takes $\jRi$ to $ \iRj$ and
$1_{\ii}$ to $1_{\ii}$.  Flipping a diagram about a vertical axis and
simultaneously taking
$$
 \xy    (0,0)*{\dcross{i}{i}}; \endxy \quad {\rm to} \quad
  -\xy (0,0)*{\dcross{i}{i}}; \endxy
$$
(in other words, multiplying the diagram by $(-1)^s$ where $s$ is the number of
times equally labelled strands intersect) is an involution $\sigma$ of $R(\nu)$
which commutes with $\psi$.

Sometimes it is convenient to convert from graphical to algebraic notation. For a
sequence $\ii=i_1 i_2\dots i_m\in \seq(\nu)$ and $1\le k \le m$ we denote
\begin{eqnarray} \label{eq_dot_xki}
  x_{k,\ii} \quad := \quad
  \xy (0,0)*{\sup{}};  (0,-6)*{\scs i_1};  \endxy
    \dots
    \xy (0,0)*{\supdot{}};  (0,-6)*{\scs i_k}; \endxy
   \dots
    \xy (0,0)*{\sup{}};  (1,-6)*{\scs i_m}; \endxy
\end{eqnarray}
with the dot positioned on the $k$-th strand counting from the left, and
\begin{eqnarray} \label{eq_crossing_delta}
  \delta_{k,\ii} \quad := \quad
  \xy  (0,0)*{\sup{}};  (0,-6)*{\scs i_1};   \endxy
    \dots
  \xy (0,0)*{\dcross{i_k}{\; \; i_{k+1}}}; \endxy
   \dots
  \xy (0,0)*{\sup{}}; (1,-6)*{\scs i_m}; \endxy
\end{eqnarray}
The symmetric group $S_m$ acts on $\seq(\nu)$, $m=|\nu|$ by permutations.
Transposition $s_k=(k,k+1)$ switches entries $i_k, i_{k+1}$ of $\ii$.  Thus,
$\delta_{k,\ii} \in {_{s_k(\ii)}R(\nu)_{\ii}}$.  The relations \eqref{eq_UUzero}
become
\begin{equation}
 \delta_{k,s_k(\ii)}\delta_{k,\ii} =
 \left\{  \begin{array}{ccc}
             0 & \text{if} & i_k=i_{k+1}. \\
             1_{\ii} & \text{if} & i_k \cdot i_{k+1}=0 \\
             x_{k,\ii}+x_{k+1,\ii} & \text{if} & i_k \cdot i_{k+1}=-1 \end{array}
 \right.
\end{equation}
Other defining relations for $R(\nu)$ can be similarly rewritten.

%
\subsection{Examples}
\label{subsec-ex}
%

\begin{enumerate}[1)]

\item $\nu =0$.  We have $R(0)=\Z$, with the unit element given by the empty
diagram.
\item $\nu=i$ for some vertex $i$.  Then a diagram is a line
with some number of dots on it.  Hence, $R(i) \cong \Z[x_{1,i}]$, where in our
notation $x_{1,i}$ denotes a line labelled $i$ with one dot on it.
\[
   \xy
  (0,0)*{\supdot{i}};
  (-10,0)*{};(10,0)*{};
 \endxy
\]

\item $\nu=m i$ for some vertex $i$. The only sequence in $\seq(mi)$ is $i^m=ii\dots
i$.  Every strand in the diagram is labelled by $i$, and the local relations are
\begin{eqnarray} \vcenter{
   \xy (0,0)*{\twocross{i}{i}}; \endxy}
  &=& 0 \qquad \qquad
 \vcenter{ \xy (0,0)*{\linecrossL{i}{i}{i}}; \endxy} =
 \vcenter{\xy  (0,0)*{\linecrossR{i}{i}{i}}; \endxy}
\end{eqnarray}
\begin{eqnarray}
\xy (0,0)*{\dcrossul{i}{i}};  \endxy
  \quad - \quad
  \xy  (0,0)*{\dcrossdr{i}{i}};  \endxy
  & = &
 \xy  (-3,0)*{\sup{i}}; (3,0)*{\sup{i}};
 \endxy \\
  \xy (0,0)*{\dcrossdl{i}{i}}; \endxy
 \quad - \quad
 \xy  (0,0)*{\dcrossur{i}{i}};  \endxy
  & = &
 \xy  (-3,0)*{\sup{i}};  (3,0)*{\sup{i}}; \endxy
\end{eqnarray}
 $R(mi)$ is isomorphic to the nilHecke ring $\BNC_m$, which is the unital ring of
endomorphisms of the abelian group $\Z[x_1, \dots , x_m]$ generated by the
endomorphisms of multiplication by $x_1, \dots, x_m$ and divided difference operators
$$ \partial_a(f(x)) = \frac{f(x)-s_a f(x)}{x_a -x_{a+1}}, \hspace{0.2in} 1\le a \le m-1,$$
where $s_a$ transposes $x_a$ and $x_{a+1}$ in the polynomial $f(x)$. The defining relations are
\[
 \begin{array}{ll}
 x_a x_b =   x_b x_a , &  \\
   \partial_a x_b = x_b\partial_a \quad \text{if $|a-b|>1$}, &
   \partial_a\partial_b = \partial_b\partial_a \quad \text{if $|a-b|>1$}, \\
  \partial_a^2 = 0,  &
   \partial_a\partial_{a+1}\partial_a = \partial_{a+1}\partial_a\partial_{a+1},  \\
   x_a \partial_a - \partial_a x_{a+1}=1,  &   \partial_a x_a - x_{a+1} \partial_a =1.
  \end{array}
\]
In the above equations $x_a$ stands for the operator of multiplication by $x_a$.
The defining relations are exactly our graphical equations on identically-colored
strands, with the relations in the first two rows corresponding to isotopies of diagrams.

The nilHecke ring is related to the theory of Schubert varieties, see~\cite{KK}, \cite{BilLak}. 
The nilCoxeter ring is the
subring of $\BNC_m$ generated by $\partial_1, \dots, \partial_{m-1}$, see
\cite[Chapter 2]{pos}, \cite{kho2}. Divided differences go back to Newton;
in the context of representation theory they appeared in~\cite{BGG}, \cite{Dem}.

The center of $\BNC_m$ is the ring of symmetric polynomials in $x_1, \dots , x_m$,
and $\BNC_m$ is isomorphic to the ring of $m!\times m!$ matrices with coefficients
in $Z(\BNC_m)$, see~\cite{Man}. Here we consider $\BNC_m$ as a graded ring,
with $\deg(\partial_a)=-2$ and $\deg(x_a)=2$. The graded nilHecke ring plays a
fundamental role in the categorification of Lusztig's quantum $\mathfrak{sl}_2$
defined by the second author~\cite{Lau1}.

For each permutation $w\in S_m$ let $\partial_w = \partial_{a_1}\dots \partial_{a_r},$
where $s_{a_1}\dots s_{a_r}$ is a minimal presentation of $w$, so that $r= l(w)$.
This element does not depend on the choice of presentation.

Define $ e_m=x_1^{m-1} x_2^{m-2}\dots x_{m-1}\partial_{w_0},$
where $w_0$ is the longest permutation. This element is an idempotent of degree $0$.
We will also use the idempotent $\psi(e_m)$ given by reflecting the diagram
of $e_m$ about the horizontal axis,
$$ \psi(e_m)= \partial_{w_0}x_1^{m-1} x_2^{m-2}\dots x_{m-1}.$$
$\BNC_m \psi(e_m)$ is a left $\BNC_m$-module isomorphic to the polynomial
representation of $\BNC_m$. The polynomial representation is the unique, up to
isomorphism and grading shifts, graded indecomposable projective $\BNC_m$-module.
The module $\BNC_m \psi(e_m)$ is nontrivial in even non-negative degrees only.

The regular representation of $\BNC_m$ decomposes as the sum of $m!$ copies
of the polynomial one. Taking the grading into account and denoting by $P_m$ the
module $\BNC_m \psi(e_m)$ with the grading shifted down by $\frac{m(m-1)}{2}$,
we get a direct sum decomposition of graded modules
$$ \BNC_m \cong P_m^{[m]!} .$$
Here $[m]!=[m][m-1]\dots [1]$ is the quantum factorial,
$[m]=\frac{{q^m}-q^{-m}}{q-q^{-1}}$,
and $M^f$ or $M^{\oplus f}$, for a graded module $M$ and a Laurent
polynomial $f=\sum f_a q^a\in \Z[q,q^{-1}]$,
denotes the direct sum over $a\in \Z$, of $f_a$ copies of $M\{a\}$.

We denote by $P_{i,m}$ the corresponding indecomposable graded projective module
over $R(mi)$ and by $e_{i,m}$ the idempotent corresponding to $e_m$ under the isomorphism
$R(mi)\cong \BNC_m$. As a graded abelian group, $P_{i,m}$ is nontrivial in degrees
$-\frac{m(m-1)}{2}+2\N$.

Letting ${}_mP$ be the right graded projective module $e_m\BNC_m\{-\frac{m(m-1)}{2}\}$,
we have a decomposition of graded right $\BNC_m$-modules
$$  \BNC_m \cong {}_mP^{[m]!} .$$
We denote by ${}_{i,m}P$ the corresponding indecomposable graded projective right
$R(mi)$-module. Idempotents $e_{i,m}$ and $\psi(e_{i,m})$ have the following diagrammatic
presentation for $m=3$
\begin{equation}
e_{i,3} \;\;=\;\;
 \xy
 (0,0)*{\linecrossR{i}{i}{i}};
 (-7.2,8)*{\bullet}; (-4.4,4.7)*{\bullet};(2.8,7.5)*{\bullet};
\endxy
\qquad \quad \psi(e_{i,3}) \;\;=\;\;
 \xy
 (0,0)*{\linecrossR{i}{i}{i}};
 (-7.2,-5)*{\bullet}; (-4.4,-1.7)*{\bullet};(1.8,-5.1)*{\bullet};
\endxy
\end{equation}
The quotient of $\Z[x_1, \dots, x_m]$ by the ideal of symmetric polynomials
is a representation $L_m$ of $\BNC_m$ which becomes irreducible upon tensoring
with any field $\Bbbk$. Over $\Bbbk$, any graded irreducible representation
of $\BNC_m$ is isomorphic to $L_m$, up to a grading shift. We denote the
corresponding irreducible representation of $R(mi)$ by $L(i^m)$.
It's nonzero in degrees $0, 2, \dots, \frac{m(m-1)}{2}$.
The representation $L(i^m)$ is isomorphic to
the representation induced from the one-dimensional graded module $L$
over $\Bbbk[x_1, \dots, x_m]$ (on which $x_1, \dots, x_m$ necessarily act
trivially).

\begin{lem} \label{lem-prekato}
The common $0$-eigenspace of operators $x_1, \dots , x_m$ on
$L(i^m) \cong \Ind (L)$ is exactly $1\otimes L$. All Jordan blocks of $x_m$
on $L(i^m)$ are of size $m$.
\end{lem}

\begin{proof} Due to the uniqueness of the irreducible module $L(i^m)$,
it is isomorphic to the module induced from the trivial representation of
the subring of $\BNC_m$ generated by the divided differences
$\partial_1, \dots, \partial_m$ and symmetric polynomials in $x_1, \dots, x_m$.
This induced representation is isomorphic, as a $\Bbbk[x_1, \dots, x_m]$-module,
to the quotient of
$\Bbbk[x_1, \dots, x_m]$ by the ideal generated by symmetric polynomials
without the constant term, and to the cohomology ring of the full flag variety.
The lemma follows from the standard facts about this quotient ring.
\end{proof}

\item $\nu=i+j$ and $i\cdot j=0$.  $\seq(i+j)=\{ij,ji\}$.  The ring $R(\nu)$ is
isomorphic to the ring of $2\times 2$ matrices with coefficients in
$\Z[x_1,x_2]$. The isomorphism is given on generators by
\[
\begin{array}{lcl}
 \xy (-3,0)*{\sup{i}}; (3,0)*{\sup{j}}; (-10,0)*{};(10,0)*{}; \endxy
 \leftrightarrow
 \left( \begin{array}{cc} 1 & 0 \\  0 & 0 \\    \end{array}  \right)
 & \qquad  &  \xy
  (-3,0)*{\sup{j}};  (3,0)*{\sup{i}};  (-10,0)*{};(10,0)*{};  \endxy
 \leftrightarrow
 \left( \begin{array}{cc} 0 & 0 \\ 0 & 1 \\    \end{array} \right)
  \\
 \xy  (0,0)*{\dcross{i}{j}};  (-10,0)*{};(10,0)*{};  \endxy
 \leftrightarrow
 \left( \begin{array}{cc} 0 & 0 \\  1 & 0 \\    \end{array} \right)
   & \qquad  &  \xy  (0,0)*{\dcross{j}{i}}; (-10,0)*{};(10,0)*{}; \endxy
\leftrightarrow
 \left( \begin{array}{cc} 0 & 1 \\  0 & 0 \\   \end{array}  \right)
  \\
 \xy  (-3,0)*{\supdot{i}};  (3,0)*{\sup{j}};  (-10,0)*{};(10,0)*{}; \endxy
 \leftrightarrow
 \left( \begin{array}{cc}  x_1 & 0 \\   0 & 0 \\   \end{array} \right)
   & \qquad  &  \xy
  (-3,0)*{\sup{i}};  (3,0)*{\supdot{j}};   (-10,0)*{};(10,0)*{};
 \endxy
 \leftrightarrow
 \left(  \begin{array}{cc} x_2 & 0 \\ 0 & 0 \\   \end{array}  \right)
\end{array}
\]
\item $\nu = i_1+ \dots +i_m$ and $i_k \cdot i_{\ell}=0$ for all $k \neq \ell$.
Then $R(\nu)$ is isomorphic to the ring of $m! \times m!$ matrices with
coefficients in $\Z[x_1, \ldots, x_m]$.  To see the isomorphism, enumerate the
rows and columns by elements of $\seq(\nu)$ and send the element $\joi$ to the
elementary $(\jj,\ii)$-matrix.

\item $\nu = \nu' + \nu''$ such that $i \cdot j =0$
for any $i \in {\rm Supp}(\nu')$ and $j \in {\rm Supp}(\nu'')$. In this case
$R(\nu)$ is isomorphic to the matrix algebra of size $\sbinom{|\nu|}{|\nu'|,
|\nu''|}$ with coefficients in $R(\nu') \otimes_{\Z} R(\nu'')$.  Indeed, except
for crossings, there are no interactions between strands from $\nu'$ and strands
from $\nu''$. A pair $\ii \in \seq(\nu')$, $\jj \in \seq(\nu'')$ defines the
sequence $\ii\jj \in \seq(\nu' + \nu'')$ and
\begin{equation}
  {_{\ii\jj}R_{\ii\jj}} \cong R(\nu') \otimes R(\nu''),
\end{equation}
since we can pull apart the $i$ and $j$ strands in any diagram $D$ with
$\ii\jj=\top(D)=\bot(D)$, using that $i_k \cdot j_{\ell}=0$, for all $k$ and
$\ell$.

\item $\nu = i +j$ and $i \cdot j= -1$.   We can identify $R(i+j)$ with the ring
of $2\times 2$ matrices with coefficients in $\Z[x_1,x_2]$ such that  the
bottom left coefficient is divisible by $x_1+x_2$:
\[
\begin{array}{lcl}
 \xy (-3,0)*{\sup{i}};  (3,0)*{\sup{j}}; (-10,0)*{};(10,0)*{}; \endxy
 \leftrightarrow
 \left( \begin{array}{cc} 1 & 0 \\ 0 & 0   \\    \end{array}   \right)
 & \qquad  &
  \xy (-3,0)*{\sup{j}}; (3,0)*{\sup{i}};  (-10,0)*{};(10,0)*{};  \endxy
 \leftrightarrow
 \left( \begin{array}{cc} 0 & 0 \\   0 & 1 \\     \end{array}  \right)
  \\
 \xy  (0,0)*{\dcross{i}{j}};  (-10,0)*{};(10,0)*{};  \endxy
 \leftrightarrow
 \left(  \begin{array}{cc}  0 & 0 \\  x_1+x_2 & 0 \\     \end{array}
 \right)
   & \qquad  &  \xy (0,0)*{\dcross{j}{i}}; (-10,0)*{};(10,0)*{};  \endxy
\leftrightarrow
 \left( \begin{array}{cc} 0 & 1 \\ 0 & 0 \\   \end{array} \right)
  \\
 \xy
  (-3,0)*{\supdot{i}};
  (3,0)*{\sup{j}};
  (-10,0)*{};(10,0)*{};
 \endxy
 \leftrightarrow
 \left(
   \begin{array}{cc}
     x_1 & 0 \\
     0 & 0 \\
   \end{array}
 \right)
   & \qquad  &  \xy
  (-3,0)*{\sup{i}};
  (3,0)*{\supdot{j}};
  (-10,0)*{};(10,0)*{};
 \endxy
 \leftrightarrow
 \left(
   \begin{array}{cc}
     x_2 & 0 \\
     0 & 0 \\
   \end{array}
 \right)
\end{array}
\]

\begin{rem}
If $i \cdot j =-1$ then the elements
\begin{eqnarray}
\xy (0,0)*{\linecrossL{i}{j}{i}}  \endxy
  \qquad {\rm and} \qquad -
\xy  (0,0)*{\linecrossR{i}{j}{i}};  \endxy
\end{eqnarray}
\end{rem}
are mutually orthogonal idempotents in $R(2i + j)$.  For instance,
\begin{eqnarray}
  \left(\; \xy  (0,0)*{\linecrossL{i}{j}{i}};  \endxy \;\right)^2
  &=&
 \vcenter{\xy
(-6,-13)*{};
 (6,8)*{};
 (0,0)*{\includegraphics[scale=0.4]{line_crossL.eps}};
 (0,15.1)*{\includegraphics[scale=0.4]{line_crossL.eps}};
 (-7,-10)*{\scs i};(0,-10)*{\scs j};(7,-10)*{\scs i};
 \endxy}
 \quad \refequal{\eqref{eq_UUzero}} \quad
 \vcenter{\xy
(-6,-13)*{};
 (6,8)*{};
 (0,0)*{\includegraphics[scale=0.4]{j1i.eps}};
 (0,15.1)*{\reflectbox{\includegraphics[angle=180, scale=0.4]{j1i.eps}}};
 (-6.3,8)*{\bullet};
 (-7,-10)*{\scs i};(0,-10)*{\scs j};(7,-10)*{\scs i};
 \endxy}
 \;\;  + \;\;
 \vcenter{ \xy
(-6,-13)*{};
 (6,8)*{};
 (0,0)*{\includegraphics[scale=0.4]{j1i.eps}};
 (0,15.1)*{\reflectbox{\includegraphics[angle=180, scale=0.4]{j1i.eps}}};
 (0.1,8)*{\bullet};
 (-7,-10)*{\scs i};(0,-10)*{\scs j};(7,-10)*{\scs i};
 \endxy} \qquad  \nn \\
&\refequal{\eqref{eq_UUzero}}&
\vcenter{ \xy
(-6,-13)*{};
 (6,8)*{};
 (0,0)*{\includegraphics[scale=0.4]{j1i.eps}};
 (0,15.1)*{\reflectbox{\includegraphics[angle=180, scale=0.4]{j1i.eps}}};
 (0.1,8)*{\bullet};
 (-7,-10)*{\scs i};(0,-10)*{\scs j};(7,-10)*{\scs i};
 \endxy}
 \quad  \refequal{\eqref{eq_iislide1}}  \quad
\vcenter{ \xy
(-6,-13)*{};
 (6,8)*{};
 (0,0)*{\includegraphics[scale=0.4]{j1i.eps}};
 (0,15.1)*{\reflectbox{\includegraphics[angle=180, scale=0.4]{j1i.eps}}};
 (5.6,-4)*{\bullet};
 (-7,-10)*{\scs i};(0,-10)*{\scs j};(7,-10)*{\scs i};
 \endxy}
 \;\; + \;\;
  \xy (0,0)*{\linecrossL{i}{j}{i}};   \endxy \nn \\
 &\refequal{\eqref{eq_UUzero}}&
 \xy  (0,0)*{\linecrossL{i}{j}{i}};  \endxy
 \label{eq_square}
\end{eqnarray}
Equation \eqref{eq_r3_hard} can be viewed as a decomposition of the idempotent
$1_{iji}$ into the sum of two orthogonal idempotents.
Equation~\eqref{eq_r3_easy} can be thought of as allowing triple intersections
for certain $ijk$.
\end{enumerate}

%
\subsection{A faithful representation and a basis of $R(\nu)$}
\label{subsection-basis}
%

{\bf Action of $R(\nu)$ on the sum of polynomial spaces.}
Choose an orientation of each edge of $\Gamma$.  For each $\nu \in \N[I] $ we
define an action of $R(\nu)$ on the free abelian group
\begin{equation}
  \Pol_{\nu} = \bigoplus_{\ii \in \seq(\nu)} \Pol_{\ii}, \qquad
  \Pol_{\ii} = \Z[x_1(\ii),x_2(\ii), \dots, x_m(\ii)], \qquad   m=|\nu| \nn.
\end{equation}
It's useful to think of the variable $x_k(\ii)$ as labelled by the vertex of
$\Gamma$ in the $k$-th position in the sequence $\ii$.  The symmetric group $S_m$
acts on $\Pol_{\nu}$ by taking $x_a(\ii)$ to $x_{w(a)}(w(\ii))$, $w \in S_m$. The
transposition $s_k$ maps $x_a(\ii)$ to $x_a(s_k(\ii))$ if $a \neq k,k+1$,
$x_k(\ii)$ to $x_{k+1}(s_k \ii)$, and $x_{k+1}(\ii)$ to $x_k(s_k \ii)$.

To define the action of $R(\nu)$, we first require that an
element $x \in \jRi$ acts by $0$ on $\Pol_{\kk}$ if $\kk \neq \ii$ and
takes $\Pol_{\ii}$ to $\Pol_{\jj}$.
We describe the action of the generators.  The dot in the $k$-th position
$x_{k,\ii}$ (see \eqref{eq_dot_xki}) acts by sending $f\in \Pol_{\ii}$ to
$x_k(\ii)f \in \Pol_{\ii}$.   The idempotent $1_{\ii}$ act by the identity on
$\Pol_{\ii}$. The crossing $\delta_{k,\ii}$ (see \eqref{eq_crossing_delta}) acts
on $f\in \Pol_{\ii}$ by
\begin{eqnarray}
  f &\mapsto & s_k f \quad \text{if $i_{k} \cdot i_{k+1}=0$}, \nn\\
  f &\mapsto & \frac{f-s_k f}{x_k(\ii)-x_{k+1}(\ii)}  \quad \text{if
  $i_{k} =i_{k+1}$}, \nn\\
  f &\mapsto & s_k f \quad \text{if} \;\;
  i_k \longleftarrow i_{k+1}, \nn \\
 f &\mapsto & (x_k(s_k\ii)+x_{k+1}(s_k\ii))(s_k f) \quad \text{if } \;\;
 i_k \longrightarrow i_{k+1}. \nn
\end{eqnarray}
Notation $i_k \longleftarrow i_{k+1}$ means that $i_{k}\cdot i_{k+1}=-1$
and this edge of $\Gamma$ is oriented from $i_{k+1}$ to $i_k$.
Note that when $i_k=i_{k+1}$ the crossing $\delta_{k,\ii}$ acts by
the divided difference operator. When all strands have the same label $i$,
the action reduces to the action of the nilHecke algebra on its polynomial
representation.

\begin{prop} \label{prop-action}
These rules define a left action of $R(\nu)$ on $\Pol_{\nu}$.
\end{prop}

\begin{proof}
We check the defining relations for $R(\nu)$. The relation  \eqref{eq_UUzero}
with $i \cdot j=0$ and $i \cdot j=-1$ and relation \eqref{eq_ijslide} are trivial to verify.
The relation~\eqref{eq_UUzero} with $i=j$ and relations \eqref{eq_iislide1},
\eqref{eq_iislide2}, and \eqref{eq_r3_easy} for $i=j=k$ are just the nilHecke
relations. The relation \eqref{eq_r3_easy} with $i,j,k$ all distinct, or  with
$i\not= j=k$,  $i\cdot j =0$, or with $i=j\not= k$,  $j\cdot k =0$ is easy to check.
The same relation with $i=j,$  $i\cdot k=-1$ or $j =k, i\cdot j=-1$ follows from the fact that
the divided difference operator annihilates symmetric polynomials.
This leaves us with the last relation \eqref{eq_r3_hard}, reproduced below
\begin{eqnarray*}
\xy (0,0)*{\linecrossL{i}{j}{i}}; \endxy
  &-&
\xy (0,0)*{\linecrossR{i}{j}{i}}; \endxy
 \quad = \quad
 \xy (-9,0)*{\up{i}}; (0,0)*{\up{j}}; (9,0)*{\up{i}}; \endxy
 \qquad \text{if $i \cdot j=-1$. }
\end{eqnarray*}
It's enough to check it on 3-stranded diagrams, with $\nu=2i+j$.
To simplify formulas, we write $x, y,$ and $z$ instead of
$x_k(\ii)$, $x_{k+1}(\ii)$, and $x_{k+2}(\ii)$
for each $\ii=\{iji, iij, jii\}$. Assume that the $ij$ edge is $i\longleftarrow j$.
The left hand side of the relation is
$$\delta_{1,jii}\delta_{2,jii}\delta_{1,iji} -
\delta_{2,iij}\delta_{1,iij}\delta_{2,iji},$$
taking $\Pol_{iji}$ to $\Pol_{iji}$. We compute the action of this element
on each monomial $x^uy^vz^w$, $u,v,w\in \N$:
 \begin{eqnarray}
& & \delta_{1,jii}\delta_{2,jii}\delta_{1,iji}(x^u y^v z^w) =
 \delta_{1,jii}\delta_{2,jii}(x^v y^u z^w) =
\delta_{1,jii} \bigg(\frac{x^vy^uz^w-x^vy^wz^u}{y-z}\bigg) \nn   \\
& & = \frac{x^uy^v z^w-x^w y^v z^u}{x-z}(x+y), \label{eq_long1}  \\
& & \delta_{2,iij}\delta_{1,iij}\delta_{2,iji}(x^u y^v z^w)
 = \delta_{2,iij}\delta_{1,iij} ((y+z)x^{u}y^{w}z^{v} ) \nn   \\
& & = \delta_{2,iij}\bigg(\frac{x^uy^{w+1}-x^{w+1}y^u}{x-y}z^v +
\frac{x^uy^w-x^w y^u}{x-y}z^{v+1}\bigg) \nn
 \\
& & = \delta_{2,iij}\bigg(\frac{x^uy^{w+1}z^v-x^{w+1}y^uz^v +
x^uy^wz^{v+1}-x^w y^uz^{v+1}}{x-y}\bigg) \nn
 \\
& & =\frac{x^uy^{v}z^{w+1}-x^{w+1}y^v z^u +
x^uy^{v+1}z^{w}-x^w y^{v+1}z^{u}}{x-z} . \label{eq_long2}
\end{eqnarray}
One can easily verify that the difference of \eqref{eq_long1} and
\eqref{eq_long2} is $x^uy^vz^w$, proving relation~\eqref{eq_r3_hard} in this case.
When $i\longrightarrow j$, we compute
 \begin{eqnarray}
& & \delta_{1,jii}\delta_{2,jii}\delta_{1,iji}(x^u y^v z^w) =
 \delta_{1,jii}\delta_{2,jii}((x+y)x^v y^u z^w)  \nn \\
& & = \delta_{1,jii} \bigg( x^{v+1}\frac{y^u z^w-y^w z^u}{y-z}+ x^v
 \frac{y^{u+1}z^w-y^w z^{u+1}}{y-z} \bigg) \nn   \\
& & = \delta_{1,jii} \bigg( \frac{x^{v+1}y^u z^w-x^{v+1}y^w z^u +
 x^v y^{u+1}z^w-x^v y^w z^{u+1}}{y-z} \bigg) \nn   \\
& & = \frac{x^{u}y^{v+1} z^w-x^{w}y^{v+1} z^u  +
 x^{u+1} y^{v}z^w-x^w y^v z^{u+1}}{x-z}, \label{eq-long3}  \\
& & \delta_{2,iij}\delta_{1,iij}\delta_{2,iji}(x^u y^v z^w)
 = \delta_{2,iij}\delta_{1,iij} (x^{u}y^{w}z^{v} )
 = \delta_{2,iij}\bigg(\frac{x^uy^{w}z^v-x^{w}y^uz^v}{x-y}\bigg) \nn \\
& & =\frac{x^uy^{v}z^w  - x^w y^{v}z^{u}}{x-z}(y+z) . \label{eq-long4}
\end{eqnarray}
Again, the difference of \eqref{eq-long3} and \eqref{eq-long4} is $x^uy^vz^w$.
Relation~\eqref{eq_r3_hard} and Proposition~\ref{prop-action} follow.
\end{proof}

{\bf A spanning set.} We look for a lower bound on the size of $R(\nu)$.
An element of this ring is a linear combination of diagrams.

If a diagram $D$ contains two strands that intersect more than once, relations
\eqref{eq_UUzero}--\eqref{eq_r3_hard} allow us to write $D$ as a linear
combination of diagrams, each with fewer intersections than $D$. Iterating, we
can write any element of $R(\nu)$ as a linear combination of diagrams
with at most one intersection between any two strands.
Furthermore, we can slide all dots in a diagram $D$ all the way to the bottom
of the diagram at the cost of adding a linear combination of diagrams with fewer
crossings than $D$. These two operations together tell us that $R(\nu)$ is
spanned by diagrams having all dots at the bottom and with each pair of strands intersecting
at most once:
\[
 \xy
 (0,0)*{\includegraphics[scale=0.5]{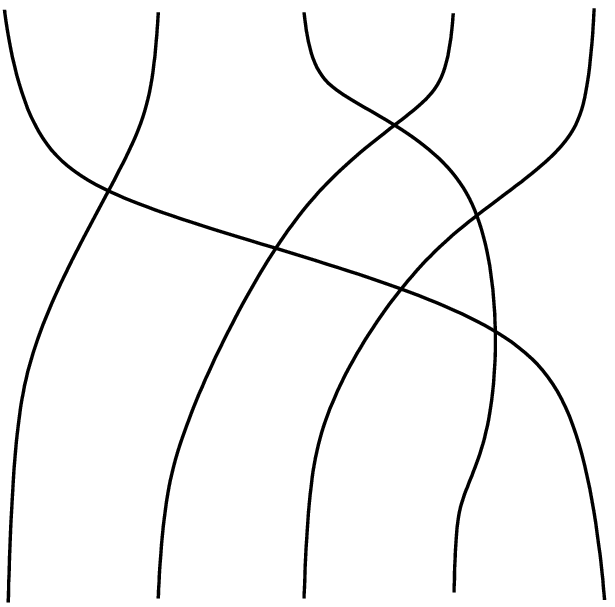}};
 (-7.2,-12)*{\bullet};(-6.7,-8)*{\bullet}; (.6,-8)*{\bullet};(14.2,-8)*{\bullet};
 (12.9,-5)*{\bullet};(14.8,-12)*{\bullet};
 \endxy
\]

Such $D$ is determined by $\ii=\bot(D)$, a minimal presentation
$\tilde{w}=s_{k_1}\dots s_{k_r}$, $r={l(w)}$
of a permutation $w\in S_m$ and the
number of dots at each bottom endpoint of $D$. The difference of two diagrams
given by the same data except for different minimal presentations of $w$ can be
written as a linear combination of diagrams with fewer crossings than each of
the original two diagrams.

For each $w\in S_m$ fix its minimal presentation $\tilde{w}$.
For $\ii, \jj\in \rm{Seq}(\nu)$ let ${_{\jj}S_{\ii}}$ be the subset of $S_m$
consisting of permutations $w$ that take $\ii$ to $\jj$ via the standard action
of permutations on sequences, defined earlier.
For each $w \in {_{\jj}S_{\ii}}$ we convert its minimal presentation $\tilde{w}$
into an element of $\jRi$ denoted $\hat{w}_{\ii}$.
Denote the subset $\{\hat{w}_{\ii}\}_{w\in {_{\jj}S_{\ii}}}$ of $\jRi$ by ${_{\jj}\hat{S}_{\ii}}$.

\begin{example} ${_{iji}S_{iji}}=\{\id,(13)\}$, and
\[
 \hat{\id}_{iji} = \xy  (-4,0)*{\up{i}}; (0,0)*{\up{j}};  (4,0)*{\up{i}}; \endxy,
 \qquad  \hat{(13)}_{iji} =
 \xy  (-4,0)*{\linecrossR{i}{j}{i}};  \endxy
 \quad {\rm or} \quad
 \xy  (-4,0)*{\linecrossL{i}{j}{i}};   \endxy  ,
\]
depending on whether we choose $s_2s_1s_2$ or $s_1s_2s_1$ as a minimal presentation
of permutation (13).
\end{example}

In general, ${_{\jj}\hat{S}_{\ii}}$ depends on our choices of minimal presentations for
permutations.  For instance, in the above example, $\hat{(13)}_{iji}\in{_{iji}\hat{S}_{iji}}$
will depend nontrivially on whether the presentation $s_1s_2s_1$ or $s_2s_1s_2$ was chosen
if $i \cdot j=-1$.

Let ${_{\jj}B_{\ii}}$ be the set $\{y \cdot x_{1,\ii}^{u_1}\dots
x_{m,\ii}^{u_m}\}$ over all $y \in {_{\jj}\hat{S}_{\ii}}$ and $u_i \in
\N$.  Here the diagrams in ${_{\jj}\hat{S}_{\ii}}$ are multiplied by all
possible monomials at the bottom. For example, the sets
${}_{ij}B_{ij}$ and ${}_{ji}B_{ij}$ consist of elements
\begin{eqnarray*}
  \xy   (-3,0)*{\supdot{i}}; (3,0)*{\supdot{j}};
  (-6,2)*{\scs u_1}; (6,2)*{\scs  u_2}; (-10,0)*{};(10,0)*{};
 \endxy
 \quad {\rm and} \qquad
  \xy   (0,0)*{\ddcross{i}{j}}; (-5.5,-1)*{\scs u_1}; (5.5,-1)*{\scs u_2};
  \endxy
\end{eqnarray*}
respectively, where we write
\begin{eqnarray*}
  \xy  (0,0)*{\supdot{}};  (-2.5,2)*{\scs u};   \endxy
      &=&
  \left(  \xy (0,0)*{\supdot{}};  \endxy  \right)^{u}.
\end{eqnarray*}

\begin{thm}\label{thm_basis}
  $\jRi$ is a free graded abelian group with a homogeneous basis ${_{\jj}B_{\ii}}$.
\end{thm}
\begin{proof}
We have already observed that the set ${_{\jj}B_{\ii}}$ spans $\jRi$. This set
consists of homogeneous elements relative to our grading on $R(\nu)$.
To prove linear independence of elements of ${_{\jj}B_{\ii}}$ we check that
they act on $\Pol_{\nu}$ by linearly independent operators.

Let ${_{\jj}1_{\ii}}$ be the diagram with the fewest number of crossings with
$\bot({_{\jj}1_{\ii}})=\ii$ and $\top({_{\jj}1_{\ii}})=\jj$.  For example,
\[
{_{jjiki}1_{ijkij}} \quad = \quad
 \xy
 (0,0)*{ \includegraphics[scale=0.4]{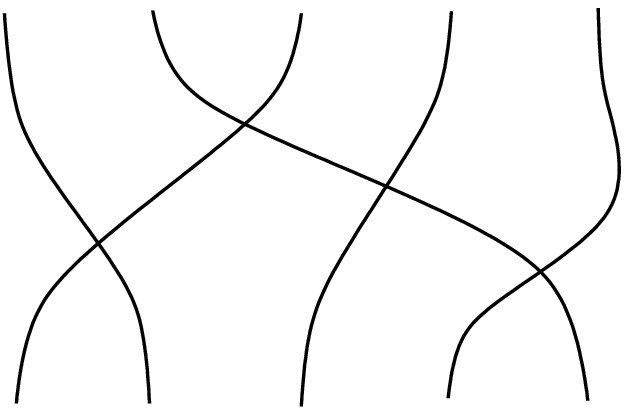}};
 (-12,-10)*{\scs i}; (-7,-10)*{\scs j};  (0,-10)*{\scs k}; (6,-10)*{\scs i}; (11.5,-10)*{\scs j};
 (-12.5,10)*{\scs j};  (-7,10)*{\scs j}; (0,10)*{\scs i}; (6,10)*{\scs k}; (12,10)*{\scs i};
 \endxy
\]
Note that identically coloured lines do not intersect in ${_{\jj}1_{\ii}}$, and that
${_{\ii}1_{\ii}}$ is just $1_{\ii}$.

The product ${_{\ii}1_{\jj}}{_{\jj}1_{\ii}} =
\prod(x_{a,\ii}+x_{b,\ii})$ where the product is over all pairs $1 \leq
a < b \leq m$ such that the lines in ${_{\jj}1_{\ii}}$ ending at $a$
and $b$ bottom endpoints counting from the left intersect and are coloured by $i,j$
with $i \cdot j =-1$.

For instance, if $\ii=ijj$, $\jj=jji$ with $i \cdot j=-1$, then
\begin{equation}
{_{\jj}1_{\ii}} \qquad = \qquad
 \xy   (0,0)*{  \includegraphics[scale=0.42]{j1i.eps}};
  (-7,-10)*{\scs i};  (0,-10)*{\scs j};  (7,-10)*{\scs j};
  (-7,10)*{\scs j};  (0,10)*{\scs j};  (7,10)*{\scs i};
 \endxy
\end{equation}
and the product
\begin{equation}
{_{\ii}1_{\jj}}{_{\jj}1_{\ii}} \quad = \quad
 \xy  (-6,-13)*{};  (6,8)*{};
   (0,0)*{\includegraphics[scale=0.4]{j1i.eps}};
   (0,15.1)*{\reflectbox{\includegraphics[angle=180, scale=0.4]{j1i.eps}}};
   (-6.5,-10)*{\scs i};(0,-10)*{\scs j};(6.5,-10)*{\scs j};
 \endxy
 \qquad = \quad (x_{1,\ii}+x_{2,\ii})(x_{1,\ii}+x_{3,\ii}).
\end{equation}

\vspace{0.1in}

Choose a complete order on the set of vertices of $\Gamma$ and orient $\Gamma$ so that
for each edge $i\longrightarrow j$
we have $i<j$ relative to the order.  This order induces a lexicographic order on
$\seq(\nu)$.  We prove linear independence of $\jBi$ by induction on $\jj\in
\seq(\nu)$ with respect to this order.

\paragraph{Base of induction:} We write
\begin{equation}
\jj = j_1 \dots j_1 j_2 \dots j_2\dots j_r\dots j_r = j_1^{\nu_1}j_2^{\nu_2}\dots
j_r^{\nu_r} \in \seq(\nu), \qquad \nu = \sum_{k=1}^r \nu_k j_k,
\end{equation}
where $j_1<j_2<\dots < j_r$ and $r$ is the cardinality of $\rm{Supp}(\nu)$.
Clearly $\jj$ is the lowest element in $\seq(\nu)$ with respect to lexicographic
order. For this $\jj$ each $w \in \jSi$ can be written uniquely as $w=w_1w_0$
where $w_1\in S_{\nu_1} \times \dots \times S_{\nu_r}$ and $w_0$ is the unique
minimal length element in $\jSi$.

Each minimal length representative $\tilde{w_0}$ determines the same
$\hat{w_0}_{\ii}=\joi$.  Likewise, the element $\hat{w_1}_{\jj}$ does not depend on the
choice of a minimal length representative $\tilde{w_1}$, since in the nilHecke
algebra the element associated to a permutation does not depend on the minimal
presentation of this permutation.
\[
 \xy
 (0,0)*{\includegraphics[scale=0.5]{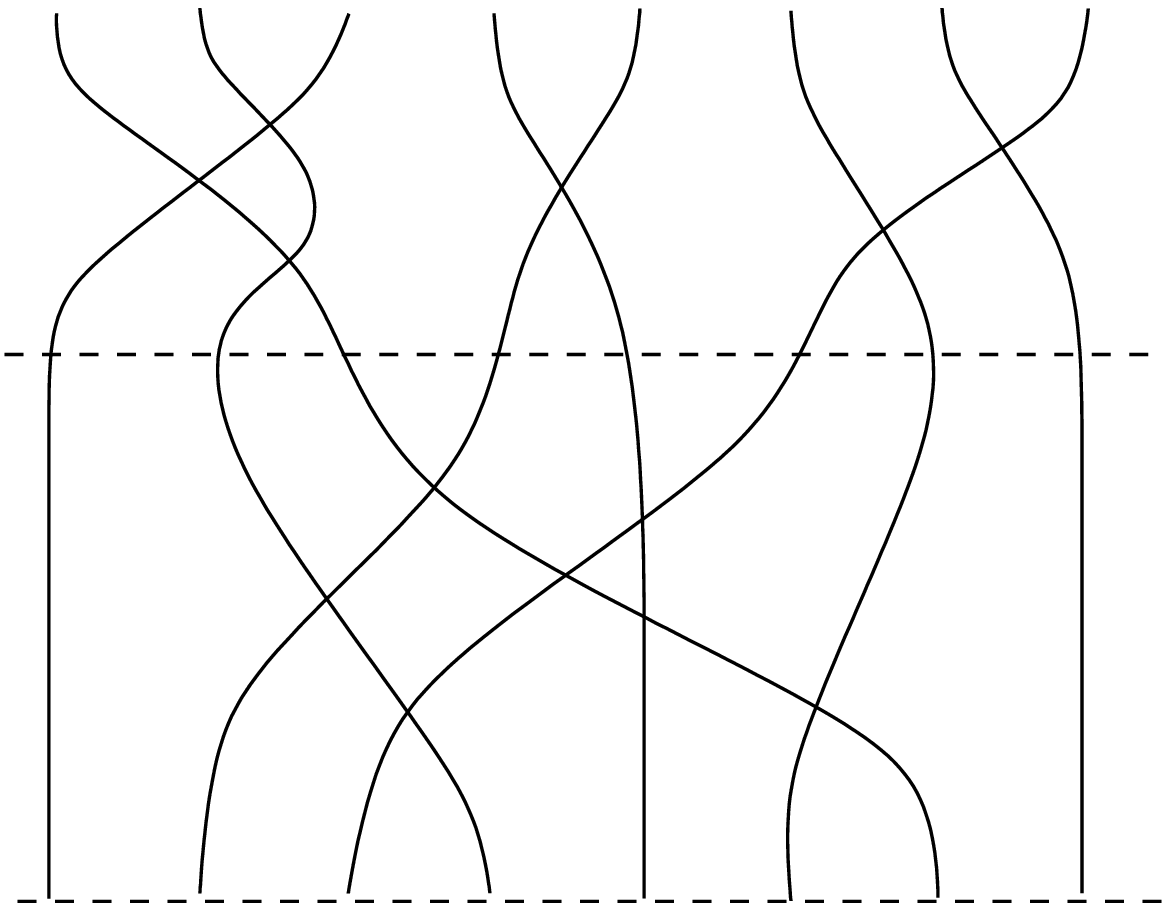}};
 (-27,25)*{\scs j_1};(-20,25)*{\scs j_1};(-12,25)*{\scs j_1};
 (-4,25)*{\scs j_2};(3,25)*{\scs j_2}; (11,25)*{\scs j_r};
 (19,25)*{\scs j_r};(26,25)*{\scs j_r};
 (35,15)*{w_1};(35,-10)*{w_0};
 \endxy
\]

The set $\jBi$ consists of elements $\hat{w_1}_{\jj}\hat{w_0}_{\ii}x^{u}$, over all
$u \in \N^m$, where
\begin{equation*}
  x^{u} = x^{u_1}_{1,\ii}x^{u_2}_{2,\ii} \dots  x^{u_m}_{m,\ii}.
\end{equation*}
It suffices to check that induced maps
\begin{equation} \label{eq_hatws}
  \hat{w_1}_{\jj}\hat{w_0}_{\ii}x^{u} \maps \Pol_{\ii} \longrightarrow \Pol_{\jj}
\end{equation}
are linearly independent, over all $u \in \N^m$.  Indeed,
$\hat{w_0}_{\ii}x^{u}$ takes $x^{v}\in \Pol_{\ii}$ to
$x^{w_0(u+v)}\in \Pol_{\jj}$, where $w_0$ acts on $v$ via the obvious
permutation. This is due to peculiarities of our action, since the
element $\delta_{k,\ii}$ takes $f \in \Pol_{\ii}$ to $s_kf$ if $i_k >i_{k+1}$.
Elements $\hat{w_1}_{\jj}$ act on the monomials by products of divided difference
operators.  It's known that the standard action of the nilHecke ring on
polynomials is faithful~\cite{Man}, implying linear independence
of all maps \eqref{eq_hatws}.

\paragraph{Induction step:} Assume we proved that $\jBi$ is independent, that
$j_k<j_{k+1}$, and set $\jj'=s_k\jj = j_1 \dots j_{k-1}j_{k+1}j_kj_{k+2}\dots
j_m$.  It suffices to assume that $j_k \cdot j_{k+1} = -1$, for otherwise $j_k
\cdot j_{k+1}=0$ and the maps $\delta_{k,\jj}$, $\delta_{k,\jj'}$ set up
bijections between $\jBi$ and ${_{\jj'}B_{\ii}}$, implying linear independence of
${_{\jj'}B_{\ii}}$.

To show that ${_{\jj'}B_{\ii}}$ is independent, we examine its image under the map
\begin{equation}
  \delta_{k,\jj'} \maps {_{\jj'}R(\nu)_{\ii}} \to \jRi.
\end{equation}
Define a partial order on $\jBi$ by requiring that $w_1x^{u}<w_2x^{v}$
if $\ell(w_1)< \ell(w_2)$, or if $w_1=w_2$, $u_1=v_1, \dots, u_t =
v_t, u_{t+1}< v_{t+1}$ for some $t$. Extend this partial order to a complete
order on $\jBi$ in some way.

Define the map $\overline{\delta} \maps {_{\jj'}B_{\ii}} \to \jBi$ by
$\overline{\delta}y=\delta_{k,\jj'}y$ if the strands of diagram $y$ ending at the top endpoints
numbered $k$, $k+1$ from the left are disjoint
\[
\overline{\delta} \;\;\maps\;\;
 \xy
 (0,0)*{\includegraphics[scale=0.5]{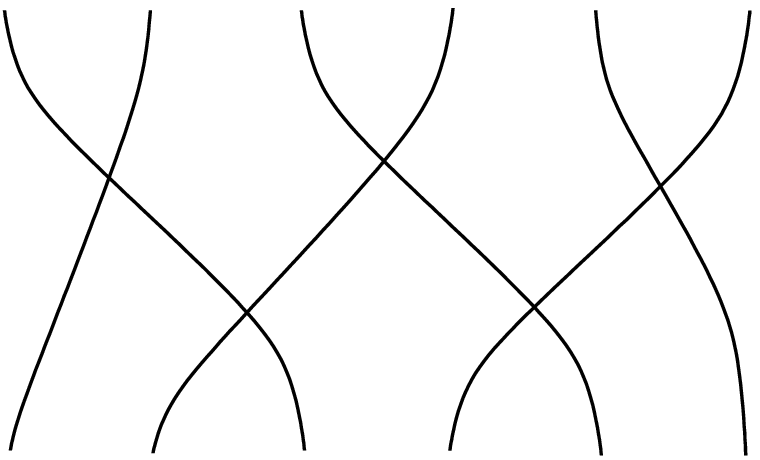}};
 (-19,15)*{\scs i_1};(4,15)*{\scs i_k};(13,15)*{\scs i_{k+1}};
 \endxy
 \quad
 \mapsto
 \quad
 \xy
 (0,0)*{\includegraphics[scale=0.5]{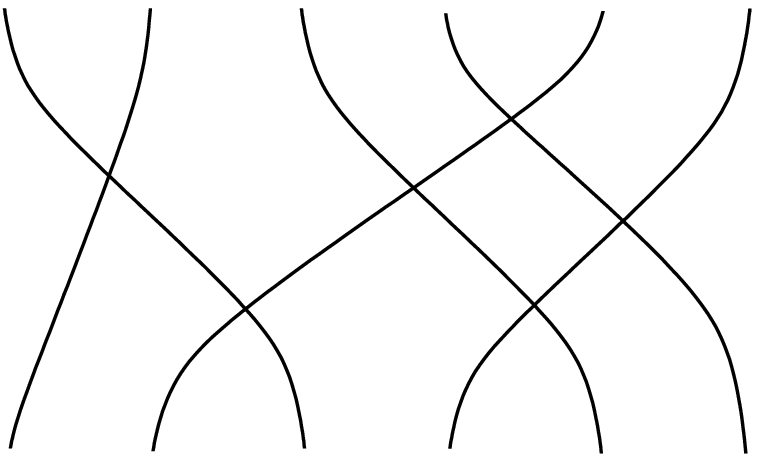}};
 (-19,15)*{\scs i_1};(4,15)*{\scs i_k};(13,15)*{\scs i_{k+1}};
 \endxy
\]
and $\overline{\delta}y = y' \cdot x_{\ell,\ii}$ if these two strands of $y$ intersect.
Here $\ell$ is the number, counting from the left, of the bottom endpoint of the strand
with top endpoint $k$, and $y'$ is obtained from $y$ by removing the intersection of
these two strands. Graphically
\[
\overline{\delta} \;\;\maps\;\;
 \xy
 (0,0)*{\includegraphics[scale=0.5]{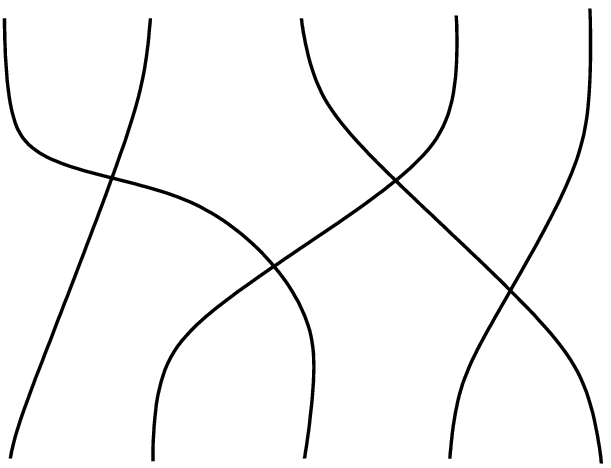}};
 (-15,14)*{\scs i_1};(0,14)*{\scs i_k};(8,14)*{\scs i_{k+1}};
 (-7.4,-8)*{\bullet}; (0.3,-8)*{\bullet};(14.5,-9)*{\bullet};
 (12.5,-5)*{\bullet};
 \endxy
 \quad
 \mapsto
 \quad
 \xy
 (0,0)*{\includegraphics[scale=0.5]{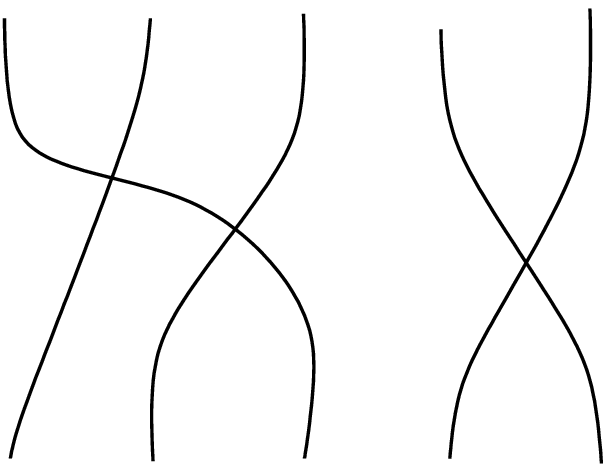}};
 (-7.8,-8)*{\bullet}; (-6,-3)*{\bullet};(0.3,-8)*{\bullet};(14.8,-9)*{\bullet};
 (13.5,-5)*{\bullet};
 (-15,14)*{\scs i_1};(0,14)*{\scs i_k};(8,14)*{\scs i_{k+1}};
 \endxy
\]

We can write $y' = \hat{s_kw}_{\ii}x^{u}x_{\ell,\ii}$ if $y=\hat{w}_{\ii}x^{u}$ as above.
The map $\overline{\delta} \maps {_{\jj'}B_{\ii}} \to \jBi$ is clearly injective.
It is not hard to compute that $\delta_{k,\jj'} y = \overline{\delta}y +$ lower terms:
\begin{equation}
  \delta_{k,\jj'} y = \overline{\delta}y +
  \sum_{z< \overline{\delta}y} n_z \cdot z, \qquad n_z \in \Z, \quad z \in {_{\jj}B_{\ii}},
\end{equation}
for any $y \in {_{\jj'}B_{\ii}}$, and by lower order terms we mean a linear
combination of elements of $\jBi$ less than $\overline{\delta}y$ with respect to
the order of $\jBi$.

Induction step follows, since $\jBi$ is a linearly independent set by induction
hypothesis.  This completes the proof of Theorem~\ref{thm_basis}.
\end{proof}

The representation $\Pol_{\nu}$ has a grading. Choose any $\ii\in \seq(\nu)$ and
place the unit element $1\in \Pol_{\ii}$ in degree $0$. This uniquely determines
a grading on $\Pol_{\nu}$ making it a graded module over the graded ring $R(\nu)$.

\begin{cor} $\Pol_{\nu}$ is a faithful graded module over the graded ring $R(\nu)$.
\end{cor}

%
\subsection{Properties of $R(\nu)$}
\label{subsection-properties}
%

From Theorem~\ref{thm_basis} we deduce several properties of $R(\nu)$.   For each
$\ii \in \seq(\nu)$ the subring $\iRi$ contains the polynomial ring
$\Pol(\nu,\ii) \cong \Z[x_{1,\ii}, x_{2,\ii}, \ldots, x_{m,\ii}]$.  We
differentiate between the ring $\Pol(\nu,\ii)$ and the abelian group $\Pol_{\ii}$
on which we defined the action.  The direct product
\begin{equation}
\Pol(\nu) = \prod_{\ii \in \seq(\nu)} \Pol(\nu,\ii)
\end{equation}
is a commutative subring of $R(\nu)$.

\begin{prop}
$R(\nu)$ is a free $\Pol(\nu)$-module of rank $m!$ with respect to both left and
right multiplication actions of $\Pol(\nu)$.
\end{prop}

\begin{proof}
For each permutation $w\in S_m$ choose a minimal representative $\tilde{w}$
and form
$$\hat{w}=\sum_{\ii\in \seq(\nu)} \hat{w}_{\ii}.$$
Theorem~\ref{thm_basis} implies that the set $\{\hat{w}\}_{w\in S_n}$
is a basis of $R(\nu)$ as a free graded module over $\Pol(\nu)$ under the right
multiplication action of the latter. The left multiplication case follows by
applying the antiinvolution $\psi$ of $R(\nu)$.
\end{proof}

The symmetric group $S_m$ acts on $\Pol(\nu)$ by permuting strands (which carry
labels and dots).  Let ${\rm Sym}(\nu) = \Pol(\nu)^{S_m}$ be the subring of
$S_m$-invariants.  It's naturally isomorphic to the tensor product of rings of
symmetric polynomials
\begin{equation}
  {\rm Sym}(\nu) \cong \bigotimes_{i \in {\rm Supp}(\nu)} \Z[x_1, \ldots,
  x_{\nu_i}]^{S_{\nu_i}},
\end{equation}
over vertices $i$ in ${\rm Supp}(\nu)$, with the number of variables $\nu_i$.

\begin{example}
  The following elements are generators of ${\rm Sym}(2i+j)$:
\[
 \begin{array}{l}
   1 =  \xy (0,0)*{\sup{i}};  (6,0)*{\sup{i}}; (12,0)*{\sup{j}} \endxy
  +
   \xy (0,0)*{\sup{i}}; (6,0)*{\sup{j}};  (12,0)*{\sup{i}} \endxy
  +
   \xy (0,0)*{\sup{j}}; (6,0)*{\sup{i}}; (12,0)*{\sup{i}} \endxy
     \\
   \xy (0,0)*{\supdot{i}{}}; (6,0)*{\sup{i}}; (12,0)*{\sup{j}} \endxy
  +
   \xy (0,0)*{\sup{i}}; (6,0)*{\supdot{i}{}}; (12,0)*{\sup{j}} \endxy
  +
   \xy (0,0)*{\supdot{i}{}}; (6,0)*{\sup{j}}; (12,0)*{\sup{i}} \endxy
  +
   \xy (0,0)*{\sup{i}}; (6,0)*{\sup{j}}; (12,0)*{\supdot{i}{}} \endxy
  +
   \xy (0,0)*{\sup{j}}; (6,0)*{\supdot{i}}; (12,0)*{\sup{i}} \endxy
  +
   \xy (0,0)*{\sup{j}}; (6,0)*{\sup{i}}; (12,0)*{\supdot{i}} \endxy
  \\
   \xy (0,0)*{\supdot{i}{}}; (6,0)*{\supdot{i}{}}; (12,0)*{\sup{j}} \endxy
  +
  \xy (0,0)*{\supdot{i}{}}; (6,0)*{\sup{j}}; (12,0)*{\supdot{i}{}} \endxy
  +
  \xy (0,0)*{\sup{j}}; (6,0)*{\supdot{i}{}}; (12,0)*{\supdot{i}{}} \endxy
  \\
   \xy (0,0)*{\sup{i}}; (6,0)*{\sup{i}}; (12,0)*{\supdot{j}{}} \endxy
  +
  \xy (0,0)*{\sup{i}}; (6,0)*{\supdot{j}{}}; (12,0)*{\sup{i}} \endxy
  +
  \xy (0,0)*{\supdot{j}{}}; (6,0)*{\sup{i}}; (12,0)*{\sup{i}} \endxy
 \end{array}
\]
\end{example}

Consider the inclusions of rings
\begin{equation}
  {\rm Sym}(\nu) \subset \Pol(\nu) \subset R(\nu) .
\end{equation}
Each subsequent ring is a free rank $m!$ module over the previous ring.
Therefore, $R(\nu)$ is a free module of rank $(m!)^2$ over ${\rm Sym}(\nu)$.
Moreover, a simple computation shows that ${\rm Sym}(\nu)$ belongs to the center
of $R(\nu)$.  The converse is true as well.

\begin{thm}
${\rm Sym}(\nu)$ is the center of $R(\nu)$.
\end{thm}

\begin{proof}
The multiplication map
\begin{eqnarray*}
  {_{\ii}R(\nu)_{\kk}} &\longrightarrow& {_{\jj}R(\nu)_{\kk}} \\
  y & \mapsto& {_{\jj}1_{\ii}y} \nn
\end{eqnarray*}
by ${_{\jj}1_{\ii}}$ is injective, since the composition
${_{\ii}1_{\jj}}{_{\jj}1_{\ii}}$ is injective, being a certain product of sums of
$x_{a,\ii}$'s (use Theorem~\ref{thm_basis}).

A central element $z \in Z(R(\nu))$ decomposes
\begin{equation}
  z = \sum_{\ii \in \seq(\nu)} z_{\ii}, \qquad z_{\ii} = 1_{\ii}z=z1_{\ii}.
\end{equation}
In particular, $z_{\ii}$ is a central element of $\iRi$.  Let $\jj =
j_1^{\nu_1}j_2^{\nu_2}\dots j_r^{\nu_r}$, for some order $j_1\ldots j_r$ of
vertices that appear in $\nu$.  The ring ${_{\jj}R(\nu)_{\jj}}$ is the tensor
product of nilHecke rings
\begin{equation}
  {_{\jj}R(\nu)_{\jj}} \cong \bigotimes_{t=1}^r \BNC_{\nu_t}
  \end{equation}
and its center is isomorphic to the tensor product of centers of nilHecke rings,
which are known to be symmetric polynomials in $\nu_t$ variables.  Moreover, the
composition
\begin{equation}
  \xy
    (-30,0)*+{{\rm Sym}(\nu)}="1";   (0,0)*+{Z(R(\nu))}="2";
    (30,0)*+{Z({_{\jj}R(\nu)_{\jj}})}="3";
    (0,-5)*+{z}="2'";  (30,-5)*+{z_{\jj}}="3'";
    {\ar "1";"2"}; {\ar "2";"3"};  {\ar@{|->} "2'";"3'"};
  \endxy
\end{equation}
is an isomorphism.

Subtracting an element of ${\rm Sym}(\nu)$, we can assume that a central
element $z$ has $z_{\jj}=0$.  Since for all $\ii$
\begin{equation}
  0= z_{\jj}({_{\jj}1_{\ii}}) = z({_{\jj}1_{\ii}}) = ({_{\jj}1_{\ii}}) z =
({_{\jj}1_{\ii}}) z_{\ii},
\end{equation}
we get $z_{\ii}=0$ since the multiplication by ${_{\jj}1_{\ii}}$ is injective.
\end{proof}

\begin{cor} \hfill
\begin{enumerate}[1)]
\item $R(\nu)$ is a free rank $(m!)^2$ module over its center ${\rm Sym}(\nu)$.
\item $R(\nu)$ is free as a graded module over ${\rm Sym}(\nu)$.
\end{enumerate}
\end{cor}

The ring ${\rm Sym}(\nu)$ is $\Z_+$-graded.  Any finitely-generated free graded
${\rm Sym}(\nu)$-module has a graded rank invariant which lies in $\N[q,q^{-1}]$.
The graded rank of the module ${\rm Sym}(\nu)\{a\}$ whose grading starts in
degree $a$ is $q^a$, and graded rank is additive under direct sum.  It's not hard
to write a combinatorial formula for the graded rank of $R(\nu)$; we leave it to
the reader as an exercise.

\begin{cor} \hfill  \begin{enumerate}[1)]
 \item $R(\nu)$ is both left and right Noetherian.
 \item $R(\nu)$ is indecomposable.
 \end{enumerate}
\end{cor}

Indecomposability is equivalent to $1$ being the only central idempotent in the
ring.  Note that $R(\nu)$ is ``almost'' positively graded. Precisely, it is zero
in degrees less that $- \sum_i\nu_i(\nu_i-1)$.

%
\subsection{Representations}
\label{subsection-reps}
%

In this and the following sections we assume that $R(\nu)$ is defined over a
field $\Bbbk$ rather than over $\Z$.
All earlier results about $R(\nu)$ remain valid over $\Bbbk$.  We view
$R(\nu)$ as a graded $\Bbbk$-algebra with every element of $\Bbbk$
in degree $0$. Let $R(\nu)\dmod$ be the category of finitely-generated graded left
$R(\nu)$-modules, $R(\nu)\fmod$ be the category of
finite-dimensional graded $R(\nu)$-modules, and $R(\nu)\pmod$ be the
category of projective objects in $R(\nu)\dmod$. The morphisms in each of
these three categories are grading-preserving module homomorphisms.
The first two categories are abelian. We have a diagram of categories and inclusions
$$ R(\nu)\fmod \subset R(\nu)\dmod \supset R(\nu)\pmod .$$
From now on, by an $R(\nu)$-module we mean a \emph{left graded finitely-generated}
$R(\nu)$-module, unless otherwise specified.
For any two $R(\nu)$-modules $M$, $N$ denote by
$\Hom(M,N)$ or $\Hom_{R(\nu)}(M,N)$ the $\Bbbk$-vector space of
grading-preserving homomorphisms, and by
\begin{equation}
  \HOM(M,N) := \bigoplus_{a \in \Z} \Hom(M,N\{a\})
\end{equation}
the $\Z$-graded $\Bbbk$-vector space of all $R(\nu)$-module morphisms.
Here $N\{a\}$ denotes $N$ with the grading shifted up by $a$.
By a simple $R(\nu)$-module we mean a simple object in the category $R(\nu)\dmod$.
We denote by ${\rm Sym}^+(\nu)$ the unique graded maximal central ideal
of ${\rm Sym}(\nu)$. It is spanned by $S_m$-invariant polynomials without the
constant term.

\begin{prop}
A simple $R(\nu)$-module $S$ is finite-dimensional and
${\rm Sym}^+(\nu)$ acts by $0$ on it. $\Hom(S,S\{a\})=0$ if $a\neq 0$, and
$S$ remains simple when viewed as an $S$-module without the grading.
\end{prop}
\begin{proof} Obvious. \end{proof}
Hence, $S$ is a (graded) module over the finite-dimensional quotient algebra
\begin{equation}  R'(\nu) = R(\nu) / {\rm Sym}^+(\nu)R(\nu). \end{equation}
Note that $\dim_{\Bbbk}R'(\nu)=(m!)^2$, and, up to isomorphism and grading
shifts, there are only finitely many simple $R(\nu)$-modules. We choose
one representative $S_b$ from each equivalence class, denote the set of
equivalence classes by $\mathbf{B}'_{\nu}$, and define
$$\mathbf{B}'\define  \bigsqcup_{\nu\in \N[I]}\mathbf{B}'_{\nu}.$$
We expect a bijection between $\mathbf{B}'$ and Lusztig-Kashiwara canonical basis $\mathbf{B}$,
hence use a similar notation. Thus, any simple
$R(\nu)$-module is isomorphic to $S_b\{a\}$ for a unique $b \in \mathbf{B}'_{\nu}$
and $a \in \Z$ (recall that we are considering only graded modules).
We do not specify the grading shift for $S_b$ yet (but see the end of
Section~\ref{subsec-injectivity}).

Each module in $R(\nu)\fmod$ has finite length composition series with subsequent
quotients---simple modules.  The Grothendieck group $G_0(R(\nu))$ of
$R(\nu)\fmod$ is a free $\Z[q,q^{-1}]$-module with the basis
$\{ [S_b]\}_{b \in \mathbf{B}'_{\nu}}$ and the multiplication by $q$ corresponding to the
grading shift up by $1$.

The abelian category $R(\nu)\dmod$  has the Krull-Schmidt unique direct sum
decomposition property for modules. Objects $P_{\ii}$, $\ii \in \seq(\nu)$, belong to its
subcategory $R(\nu)\pmod$ of projective modules.

Each simple $S_b$ has a unique (up to isomorphism) indecomposable projective cover,
denoted $P_b$. We have $\HOM(P_b,S_b) \cong \Hom(P_b,S_b)\cong\End(S_b).$
An indecomposable object
of $R(\nu)\pmod$ is isomorphic to $P_b\{a\}$ for a unique $b \in \mathbf{B}'_{\nu}$
and $a\in \Z$. Any object of $R(\nu)\pmod$ has a unique, up to isomorphism, direct
sum decomposition into indecomposables.  The Grothendieck group of $R(\nu)\pmod$
is a free $\Z[q,q^{-1}]$-module with the basis $\{ [P_b]\}_{b\in \mathbf{B}'_{\nu}}$
given by the images $[P_b]$ of indecomposable projectives. Denote this Grothendieck
group by $K_0(R(\nu))$.

Recall that for a right, respectively left, $R(\nu)$-module $M$ we denote by
$M^{\psi}$ the left, respectively right $R(\nu)$-module $M$ with the action twisted
by $\psi$. For $P\in R(\nu)\pmod$, let $\overline{P}= \HOM(P, R(\nu))^{\psi}$.
This is a graded projective left $R(\nu)$-module and $\bar{\ }$
is a contravariant self-equivalence in $R(\nu)\pmod$. We have
$\overline{P_{\ii}} \cong P_{\ii}$ for each $\ii\in\seq(\nu)$, and, more generally,
$\overline{P_{\ii}\{a\}} \cong P_{\ii}\{-a\}.$ This self-equivalence induces
a $\Z[q,q^{-1}]$-antilinear involution on $K_0(R(\nu))$, also denoted $\bar{\ \ }$.

\vspace{0.1in}

There is a $\Z[q,q^{-1}]$-bilinear pairing
\begin{eqnarray} \label{eq_bil_pair}
& &  (,) \maps K_0(R(\nu)) \times G_0(R(\nu)) \longrightarrow \Z[q,q^{-1}], \\
 & &  ([P],[M]) := {\rm gdim}_{\Bbbk}(P^{\psi}\otimes_{R(\nu)}M).
\end{eqnarray}
When the field $\Bbbk$ is algebraically closed, $\End(S_b) \cong \Bbbk$,
and the bases $\{ [P_{b}]\}_b$ and $\{[S_b]\}_b$ are dual, possibly up to
rescaling by powers of $q$ and permutation of elements.
In this case $G_0(R(\nu))$ and $K_0(R(\nu))$ are dual free $\Z[q,q^{-1}]$-modules.
We will show in Section~\ref{subsec-injectivity} that, over any field $\Bbbk$,
simples $S_b$ are absolutely irreducible and the above pairing is perfect without
any restrictions on $  \Bbbk$.

There is a $\Z[q,q^{-1}]$-bilinear form, also denoted $(,)$,
 \begin{equation} \label{eq_bil_pair2}
  (,)\maps K_0(R(\nu)) \times K_0(R(\nu)) \longrightarrow \Z[q^{-1},q]\cdot (\nu)_q,
\end{equation}
where
\begin{equation}
  (\nu)_{q} = {\rm gdim} ({\rm Sym}(\nu)) = \prod_{i \in \Gamma}
  \left(  \prod_{a=1}^{\nu_i} \frac{1}{1-q^{2a}} \right) ,
\end{equation}
and
\begin{equation} \label{eq-bil-ten}
([P],[Q]) = {\rm gdim}_{\Bbbk}(P^{\psi}\otimes_{R(\nu)}Q ).
\end{equation}
Since $P^{\psi}\otimes_{R(\nu)}Q\cong Q^{\psi}\otimes_{R(\nu)} P$,
the form is symmetric.
It follows from Theorem~\ref{thm_basis} that
$\iRj\cong {}_{\ii}P\otimes_{R(\nu)} P_{\jj}$ is a free graded ${\rm
Sym}(\nu)$-module for any $\ii$, $\jj$.  Therefore, $P^{\psi}\otimes_{R(\nu)}Q$
is a free graded ${\rm Sym}(\nu)$-module of finite rank for any $P$, $Q$ as above,
and the form takes values in $\Z[q^{-1},q]\cdot (\nu)_q$.
We have
$$ ([P_{\jj}], [P_{\ii}]) = \gdim({}_{\jj} P \otimes_{R(\nu)} P_{\ii}) =
\gdim ({}_{\jj}R(\nu)_{\ii}).$$

\vspace{0.1in}

Define the character $\chr(M)$ of a graded finitely-generated $R(\nu)$-module $M$
as
$$ \chr(M) = \sum_{\ii\in \seq(\nu)} \gdim (1_{\ii} M) \cdot \ii.$$
The character is an element of the free $\Z((q))$-module with the basis $\seq(\nu)$;
when $M$ is finite-dimensional, $\chr(M)$ is an element of the free $\Z[q,q^{-1}]$-module
with basis $\seq(\nu)$.
We abbreviate $\gdim(1_{\ii}M)$ to $\chr(M,\ii)$,
$$ \chr(M) = \sum_{\ii\in \seq(\nu)} \chr(M,\ii)\cdot \ii .$$

Let $\pseq(\nu)$ be the set of all expressions $i_1^{(n_1)}i_2^{(n_2)}\dots i_r^{(n_r)}$
such that $n_1, \dots, n_r\in \N$ and $\sum_{a=1}^r n_a i_a=\nu$.
For instance, $$ \pseq(2i+j) = \{iij, iji, jij, i^{(2)} j, ji^{(2)}\}.$$
To $\ii\in \pseq(\nu)$ we assign the idempotent
$$1_{\ii}  = e_{i_1, n_1}\otimes e_{i_2, n_2}\otimes \dots \otimes e_{i_r, n_r},$$
given by the tensor product of minimal idempotents
$e_{i,n}$  in the nilHecke rings, see Section~\ref{subsec-ex}.

Let $\ii ! = [n_1]! \dots [n_r]! $, and
 $\hat{\ii}$ be the element of $\seq(\nu)$ given by expanding $\ii$,
$$ \hat{\ii} = i_1\dots i_1 i_2\dots i_2 \dots i_r \dots i_r.$$
$\hat{\ii}=\ii$ iff $\ii\in \seq(\nu).$
We have the equality of graded dimensions
$$ \gdim(1_{\hat{\ii}}M) = q^{-\langle\ii\rangle}\ii ! \cdot   \gdim (1_{\ii}M), \ \
\langle \ii \rangle = \sum_{k=1}^r \frac{n_k(n_k-1)}{2} ,
 $$
which follows from the structure of the nilHecke algebra. Let
$$\chr(M, \ii) = q^{ -\langle \ii \rangle}\cdot\gdim (1_{\ii}M), $$
then
$$ \chr(M, \hat{\ii}) = \ii ! \cdot \chr(M, \ii) .$$
In particular, $\chr(M)$ determines $\chr(M,\ii)$ for any $\ii \in \pseq(\nu)$.

For $\ii\in \pseq(\nu)$ define the left graded projective module
$$P_{\ii}  = R(\nu) \psi(1_{\ii}) \{ - \langle \ii \rangle \}, $$
and the right graded projective module
$$ {}_{\ii}P = 1_{\ii}R(\nu) \{- \langle \ii \rangle\} .$$
 We have $ P_{\hat{\ii}} \cong P_{\ii}^{\ii !}$ and
$ {}_{\hat{\ii}}P \cong {}_{\ii}P^{\ii !}$ .  For instance,
$$P_{ii}\cong P_{i^{(2)}}^{[2]!}= P_{i^{(2)}}^{q+q^{-1}}=
P_{i^{(2)}} \{1\} \oplus P_{i^{(2)}}\{ -1\}.$$
Moreover,
\begin{equation} \label{eq_char1}
\chr(M, \ii) = \gdim ({}_{\ii}P\otimes_{R(\nu)}M) = \gdim( \rm{HOM}(P_{\ii}, M)).
\end{equation}

Given two or more sequences in $\pseq(\nu)$ that differ only in several neighbouring terms,
we denote identical parts in them via dots. For instance,
$ \dots i j \dots $ and $\dots ji \dots$ denote a pair of sequences
${\ii'}ij {\ii''}$ and ${\ii'}ji{\ii''}$ for some sequences $\ii'$, $\ii''$.

\begin{prop} \label{prop-iso-leftright}
There are isomorphisms of graded projective right $R(\nu)$-modules
\begin{eqnarray*}
 {}_{\dots ij \dots}P & \cong & {}_{\dots, ji\dots}P  \ \ \textrm{if} \ \
 i\cdot j = 0, \\
 {}_{\dots iji\dots}P  & \cong & P_{\dots i^{(2)} j \dots}\oplus
 P_{\dots j i^{(2)} \dots} \  \  \textrm{if} \  \  i \cdot j = -1 ,
\end{eqnarray*}
and isomorphisms of graded projective left $R(\nu)$-modules
\begin{eqnarray*}
 P_{\dots ij \dots} & \cong & P_{\dots ji\dots} , \ \ \textrm{if} \ \
 i\cdot j = 0, \\
 P_{\dots iji\dots} & \cong & P_{\dots i^{(2)} j \dots}\oplus
 P_{\dots j i^{(2)} \dots} \  \  \textrm{if} \  \  i \cdot j = -1 .
\end{eqnarray*}
\end{prop}

\begin{proof} It suffices to show the isomorphisms for right projective modules;
application of the antiinvolution $\psi$ would imply the corresponding isomorphisms
for left projective modules. Multiplication by the diagram
\begin{eqnarray*}
  \xy  (0,0)*{\sup{}};  (0,-6)*{\scs };   \endxy
    \dots
  \xy (0,0)*{\dcross{i}{\; \; j}}; \endxy
   \dots
  \xy (0,0)*{\sup{}}; (1,-6)*{\scs }; \endxy
\end{eqnarray*}
is a grading-preserving  isomorphism between ${}_{\dots ij\dots} P$ and ${}_{\dots ji \dots}P$
if $i\cdot j=0$.

Consider grading-preserving maps
\begin{eqnarray*}
B_0 & : & {}_{\dots iji \dots} P \lra {}_{\dots i^{(2)} j \dots}P\oplus
   {}_{\dots j i^{(2)} \dots}P   \\
B_1 & : & {}_{\dots i^{(2)} j \dots}P \oplus
   {}_{\dots j i^{(2)} \dots}P \lra {}_{\dots iji \dots} P
\end{eqnarray*}
given by matrices of diagrams
\begin{equation}
\eqdefB
\end{equation}

Notice that the top entry in $B_0$ ends with the projector $e_{2,i}$ which takes
 $1_{\dots ii \dots }$ to $1_{\dots i^{(2)} \dots}$, and we view this entry as a
homomorphism
$${}_{\dots iji \dots} P \lra {}_{\dots i i j \dots}P \lra {}_{\dots i^{(2)} j \dots}P,$$
ditto for the bottom entry in $B_0$. The degree of each diagram in $B_0$ is $1$,
therefore the map $B_0$ is grading-preserving since the grading shifts for the
sequences differ by
$ \langle {\dots i^{(2)} j \dots} \rangle - \langle  {\dots iji \dots} \rangle= 1.$

We view the first entry in $B_1$ as the composition
$$ {}_{\dots i^{(2)} j \dots}P \subset {}_{\dots i i j \dots} P \lra {}_{\dots i j i\dots} P,$$
and there is no need to write the corresponding idempotent (same for the second entry).
We compute
\begin{eqnarray} \BoBoneA \nn \end{eqnarray}
\begin{eqnarray}\BoBoneB \nn \end{eqnarray}

In the last matrix above the diagonal terms are idempotents $1_{i^{(2)}j}$ and
$1_{j i^{(2)}}$, giving identity maps of projectives
${}_{\dots i^{(2)} j \dots}P$ and ${}_{\dots j i^{(2)}\dots}P$, respectively.

\begin{eqnarray} \BoneBoA  \nn \end{eqnarray}
\begin{eqnarray}\BoneBoB \nn \end{eqnarray}

Therefore, $B_0, B_1$ are isomorphisms, and the second isomorphism  follows.
\end{proof}

\begin{cor} For any $M$ in $R(\nu)\dmod$ there are isomorphisms of
graded vector spaces
\begin{eqnarray*}
 1_{\dots ij \dots}M & \cong & 1_{\dots ji\dots} M \ \ \textrm{if} \ \
 i\cdot j = 0, \\
 1_{\dots iji\dots}M\{1\} & \cong & 1_{\dots i^{(2)} j \dots}M \oplus
1_{\dots j i^{(2)} \dots}M \  \  \textrm{if} \  \  i \cdot j = -1 .
\end{eqnarray*}
\end{cor}

\begin{cor} The following character equalities hold for
any graded finitely-generated $R(\nu)$-module $M$
\begin{eqnarray*}
& & \chr(M,\dots ij \dots ) = \chr(M, \dots, ji\dots)  \ \ \textrm{if} \ \
 i\cdot j = 0 , \\
& & \chr(M, \dots iji\dots )= \chr (M, \dots i^{(2)} j \dots ) +
\chr(M, \dots j i^{(2)} \dots) \ \textrm{if}  \ i \cdot j = -1 , \\
& & \chr(M, \dots i^{(a)} i^{(b)} \dots ) =
\left[ \begin{array}{c} a+b \\ a  \end{array}
  \right] \chr(M, \dots i^{(a+b)} \dots ).
\end{eqnarray*}
\end{cor}

%
\subsection{Induction and restriction}
\label{subsection-indres}
%

Suppose we have an inclusion of rings $\iota: B \hookrightarrow A$ which is
not necessarily unital: $e=\iota(1)$ is only an idempotent in $A$. The induction
functor between categories of \emph{unital} modules
$$ B\dmod \stackrel{\Ind}{\lra}
 A\dmod , \quad M\longmapsto A\otimes_B M $$
is isomorphic to the functor $M\longmapsto Ae\otimes_B M.$
Its right adjoint
$$\Res : A\dmod \lra B\dmod $$
takes $M$ to $eM$, viewed as a $B$-module.

The inclusion of graded rings
$$\iota_{\nu,\nu'}\ : \  R(\nu)\otimes R(\nu') \hookrightarrow R(\nu+\nu')$$
is described by putting the diagrams next to each other. It takes the
idempotent $1_{\ii}\otimes 1_{\jj}$ to $1_{\ii\jj}$ and the unit
element to an idempotent of $R(\nu+\nu')$ denoted $1_{\nu,\nu'}$.

\begin{prop}  $1_{\nu,\nu'}R(\nu+\nu')$ is
a free graded left $R(\nu)\otimes R(\nu')$-module.
\end{prop}

\begin{proof} The minimal representative $w$ of a left $S_{|\nu|}\times S_{|\nu'|}$-coset
in $S_{|\nu|+|\nu'|}$ gives rise to the diagram
$${}_{\ii\jj}\hat{w}\in {}_{\ii\jj}R(\nu +\nu')_{w^{-1}(\ii\jj)}$$
of the minimal presentation of $w$ with top ends of strands labelled by the sequence
$\ii\jj$ for $\ii\in\seq(\nu)$ and $\jj\in\seq(\nu')$. 
\[
 \xy
 (0,0)*{\includegraphics[scale=0.5]{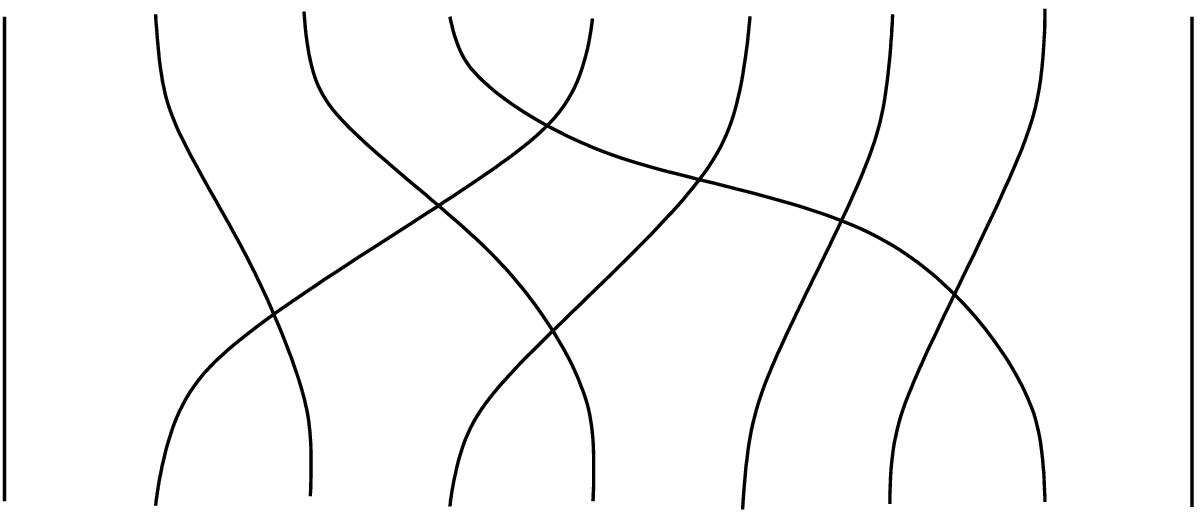}};
 (-30,16)*{ i_1};(-22,16)*{ i_2};(-8,16)*{ i_n};
 (0,16)*{ j_1};(7,16)*{ j_2};(32,16)*{ j_{m-n}};
 \endxy
\]

The set of elements
 $$\hat{w}= \sum_{\ii\in \nu, \jj\in \nu'} {}_{\ii\jj}\hat{w},$$
over all cosets, is a basis of $1_{\nu,\nu'}R(\nu+\nu')$ as a free graded left
$R(\nu)\otimes R(\nu')$-module.
\end{proof}

We denote the restriction and induction functors for the inclusion
$$  R(\nu)\otimes R(\nu') \subset R(\nu+\nu')$$
by $\Res_{\nu,\nu'}$ and $\Ind_{\nu,\nu'}$, respectively.

\begin{cor} \label{cor-res-proj} The restriction functor
$\Res_{\nu,\nu'}$ takes projectives to projectives.
\end{cor}

Given a quadruple $(\nu,\nu',\nu'',\nu''')$ with $\nu+\nu'=\nu''+\nu'''$,
let
$${}_{\nu,\nu'}R_{\nu'',\nu'''}=1_{\nu}\otimes 1_{\nu'} R(\nu+\nu')1_{\nu''}\otimes
 1_{\nu'''}.$$

\begin{prop} \label{prop-bimod-filtration}
Graded $(R(\nu)\otimes R(\nu'), R(\nu'')\otimes R(\nu'''))$-bimodule
${}_{\nu,\nu'}R_{\nu'',\nu'''}$ has a filtration by graded bimodules isomorphic to
$$  \big( {}_{\nu}R_{\nu-\lambda,\lambda} \otimes
  {}_{\nu'}R_{\nu'+\lambda-\nu''',\nu'''-\lambda} \big)
 \otimes_{R'}
 \big( {}_{\nu-\lambda,\nu''+\lambda-\nu}R_{\nu''} \otimes
  {}_{\lambda, \nu'''-\lambda}R_{\nu'''}\big) \{-\lambda\cdot (\nu'+\lambda-\nu''')\}, $$
where $R'=R(\nu-\lambda)\otimes R(\lambda)\otimes R(\nu'+\lambda-\nu''')\otimes
R(\nu'''-\lambda)$,
over all $\lambda\in \N[I]$ such that every term above is in $\N[I]$.
\end{prop}

\begin{proof} This proposition is a version of the Mackey's induction-restriction
theorem for inclusion of maximal parabolic subgroups $S_{m-n}\times S_n\subset S_m$.
The statement and its proof are best illustrated by the diagram
\[
 \xy
 (0,0)*{\includegraphics[scale=0.5]{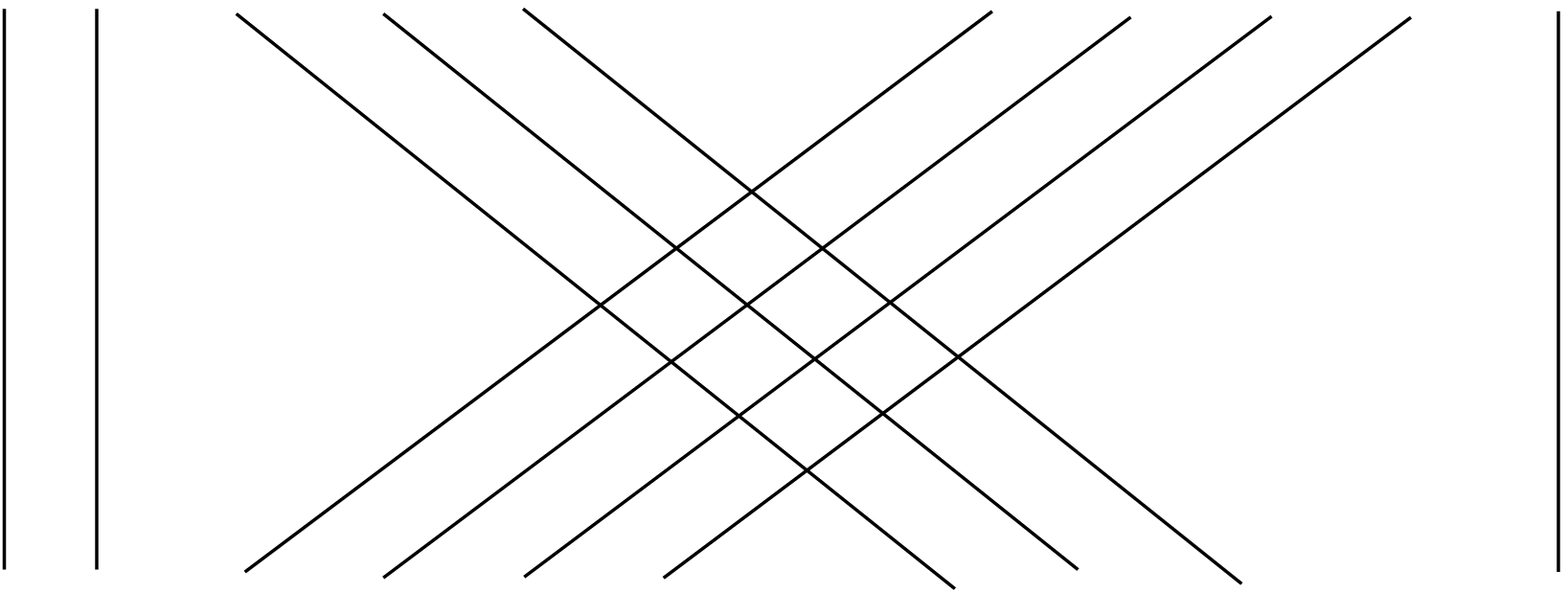}};
  (-24,-18)*{\underbrace{\hspace{1.6in}}};
  (24,-18)*{\underbrace{\hspace{1.5in}}};
  (-24,18)*{\overbrace{\hspace{1.6in}}};
  (24,18)*{\overbrace{\hspace{1.5in}}};
  (-24,22)*{\nu};(25,22)*{\nu'};(-24,-22)*{\nu''};(25,-22)*{\nu'''};
  (-50,-5)*{\nu-\lambda};(50,-5)*{\nu'''-\lambda};(-28,10)*{\lambda};(26,0)*{\nu'+\lambda-\nu'''};
 \endxy
\]

These diagrams, over all $\lambda$ (summing over all colorings of strands) will provide
generators for the subquotient bimodules that appear in the proposition. The grading
shift $-\lambda\cdot (\nu'+\lambda-\nu''')$ is the degree of the intersection diagram of
$|\lambda|$ parallel lines colored by any $\ii\in \seq(\lambda)$ and
$|\nu'+\lambda-\nu'''|$ parallel lines colored by any $\jj\in \seq(\nu'+\lambda-\nu''')$.
\end{proof}

We have
$$\Ind_{\nu,\nu'} (P_{\ii}\otimes P_{\jj}) \cong P_{\ii\jj}$$
for $\ii\in \seq(\nu), \jj\in \seq(\nu').$ By passing to direct summands,
we see that the formula holds more generally, for $\ii\in \pseq(\nu), \jj\in \pseq(\nu').$

A shuffle $\kk$ of a pair of sequences $\ii\in \seq(\nu), \jj\in \seq(\nu')$ is a
sequence together with a choice of subsequence isomorphic to $\ii$ such that $\jj$ is the
complementary subsequence. Shuffles of $\ii,\jj$ are in a bijection with the minimal
coset representatives of $S_{|\nu|}\times S_{|\nu'|}$ in $S_{|\nu|+|\nu'|}$.
We denote by $\deg(\ii,\jj,\kk)$ the degree of the diagram in $R(\nu+\nu')$
naturally associated to the shuffle, see an example below.
\[
 \xy
 (0,0)*{\includegraphics[scale=0.5]{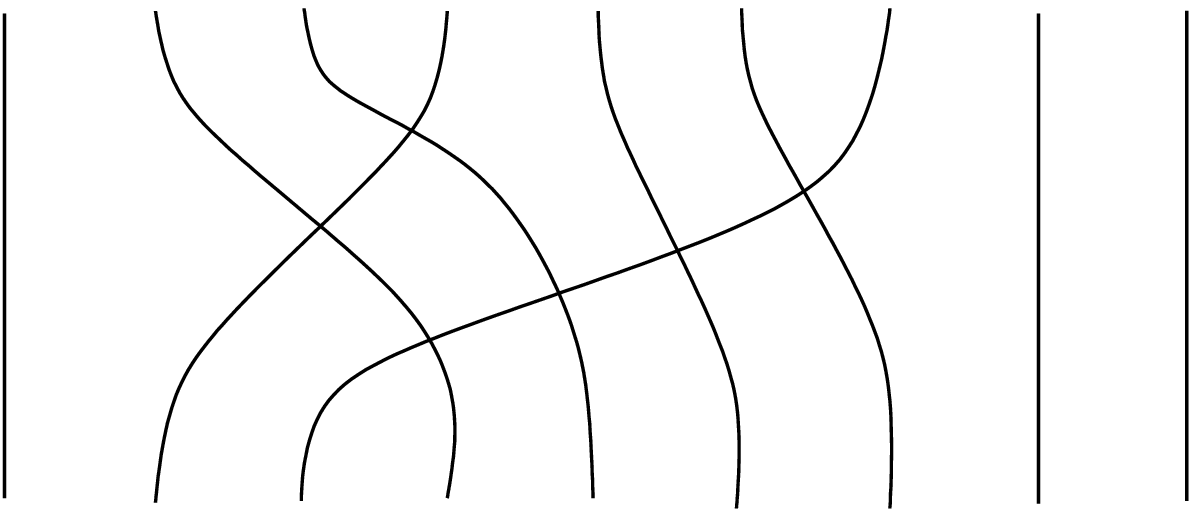}};
 (0,17)*{\overbrace{\hspace{2.5in}}};
 (-23,-17)*{\underbrace{\hspace{.75in}}};
 (11,-17)*{\underbrace{\hspace{1.6in}}};
 (0,21)*{\kk};(-23,-21)*{\ii};(11,-21)*{\jj};
 \endxy
\]
When the meaning is clear, we will also denote by $\kk$ the underlying sequence
of the shuffle $\kk$.

\begin{prop} \label{prop_res_proj}
For any $\kk\in \seq(\nu+\nu')$
\begin{eqnarray*}
\Res_{\nu,\nu'} P_{\kk} & \cong & \bigoplus_{\ii,\jj} P_{\ii}\otimes P_{\jj}\{\deg(\ii,\jj,\kk)\}, \\
\Res_{\nu,\nu'} (_{\kk}P) & \cong&  \bigoplus_{\ii * \jj = \kk}  {}_{\ii}P\otimes {}_{\jj} P
 \{\deg(\ii,\jj,\kk)\},
\end{eqnarray*}
the sum over all ways to represent $\kk$ as a shuffle of $\ii\in \seq(\nu)$ and
$\jj \in \seq(\nu')$.
\end{prop}
The proposition follows immediately from the structure of 
bimodules ${}_{\nu,\nu'}R_{\nu+\nu'}$ and
${}_{\nu+\nu'}R_{\nu,\nu'}$. $\square$

Given two functions $f$ and $g$ on sets $\seq(\nu)$ and $\seq(\nu')$, respectively,
with values in some commutative ring which contains $\Z[q,q^{-1}]$, we
define their (quantum) shuffle
product $f\shuffle g$ (see~\cite{Lec} and references therein) 
as a function on $\seq(\nu+\nu')$ given by
$$ (f\shuffle g)(\kk) = \sum_{\ii,\jj} q^{\deg(\ii,\jj,\kk)} f(\ii) g(\jj),$$
the sum is over all ways to represent $\kk$ as a shuffle of $\ii$ and $\jj$.

\begin{lem} \label{lemma_shuffle}
For $M\in R(\nu)\dmod$ and $N\in R(\nu')\dmod$ we have
$$\chr(\Ind_{\nu,\nu'}(M\otimes N) ) = \chr(M) \shuffle \chr(N).$$
\end{lem}
\begin{proof}
This lemma follows at once from the last proposition and formula~\eqref{eq_char1}. 
\end{proof} 

%
%
\section{Quantum groups and the Grothendieck ring of $R$}
\label{sec-algebras}
%

%
\subsection{Homomorphism $\gamma$ of twisted bialgebras}
\label{subsec-dir-sum}
%

For a graph $\Gamma$, we form the direct sum
$$ R = \bigoplus_{\nu\in \N[I]} R(\nu). $$
This is a non-unital ring. By various categories of $R$-modules we will
mean direct sums of corresponding categories of $R(\nu)$-modules:
\begin{eqnarray*}
R\dmod & \define &  \bigoplus_{\nu\in \N[I]} R(\nu)\dmod , \\
R\fmod & \define & \bigoplus_{\nu\in \N[I]} R(\nu)\fmod , \\
R\pmod & \define & \bigoplus_{\nu\in \N[I]} R(\nu)\pmod.
\end{eqnarray*}
The Grothendieck groups
$$ K_0(R) = \bigoplus_{\nu\in \N[I]} K_0(R(\nu)), \ \ \
     G_0(R) = \bigoplus_{\nu\in \N[I]} G_0(R(\nu)) $$
are the direct sums of Grothendieck groups of rings $R(\nu)$.
We extend the pairings~\eqref{eq_bil_pair} and~\eqref{eq_bil_pair2}
to $K_0(R)$ and $G_0(R)$ by
requiring that subspaces corresponding to different $\nu$'s be orthogonal.
Induction and restriction functors
for the inclusion $R(\nu)\otimes R(\nu')\subset R(\nu + \nu')$,
summed over all $\nu, \nu'$, give  functors
$$\Ind: R\otimes R\dmod \lra R\dmod, \ \ \  \Res:
  R\dmod \lra R\otimes R\dmod $$
where by $R\otimes R\dmod$ we mean the direct sum of categories
$R(\nu)\otimes R(\nu')\dmod,$ over all $\nu, \nu'$. These functors
restrict to subcategories of finite-dimensional modules and projective
modules. Indeed, induction takes
projectives to projectives. Restriction, in the case of these inclusions,
also takes projectives to projectives, by Proposition~\ref{prop_res_proj}
and the Krull-Shmidt property.

Thus, these functors induce maps $[\Ind], [\Res]$ on Grothendieck groups
$K_0(R)$ and $G_0(R)$. Note that $[\Res]$ is the sum of maps
$$ K_0(R(\nu+ \nu')) \lra K_0(R(\nu)) \otimes K_0(R(\nu')),$$
or
$$ G_0(R(\nu+ \nu')) \lra G_0(R(\nu)) \otimes G_0(R(\nu')),$$
over all $\nu, \nu'$; the tensor products here and further are over
$\Z[q,q^{-1}]$.

\begin{prop} $[\Ind]$ turns $K_0(R)$ and $G_0(R)$ into
associative unital $\Z[q,q^{-1}]$-algebras. $[\Res]$ turns
$K_0(R)$ and $G_0(R)$ into coassociative counital $\Z[q,q^{-1}]$-coalgebras.
\end{prop}
 Proof follows from the associativity of induction and restriction. The unit
element is given by inducing with the one-dimensional module over $R(\emptyset)$.
The counit is given by restricting to $R(\emptyset)$ and taking the graded
dimension. $\square$

Denote the product $[\Ind](x_1,x_2)$ for $x_1,x_2\in K_0(R)$ simply by $x_1x_2$.

We equip $K_0(R)\otimes K_0(R)$ with the algebra structure via
\begin{eqnarray} \label{eq_quasi_commute}
(x_1\otimes x_2) (x_1'\otimes x_2') = q^{ -|x_2|\cdot |x_1'|} x_1 x_1' \otimes
x_2 x_2'
\end{eqnarray}
for homogeneous $x_1, x_2, x_1', x_2'$, where $|x_2|\in \N[I]$ is the
weight of $x_2$, etc.

\begin{prop} $[\Res]$ is an algebra homomorphism from
$K_0(R)$ to $K_0(R)\otimes K_0(R)$ with the above algebra structure.
\end{prop}
\begin{proof} 
This follows from Proposition~\ref{prop-bimod-filtration}. 
\end{proof}

Recall the symmetric  bilinear pairing~\eqref{eq_bil_pair2} on $K_0(R)$
taking values in $\Z[q,q^{-1}]\cdot (\nu)_q.$
\begin{prop} \label{prop-pairing-prop}
The pairing $(,)$ has the following properties
\begin{enumerate}
 \item $(1,1)=1$,
\item $([P_i], [P_j]) = \delta_{i,j}(1-q^2)^{-1} $ for $i,j\in I$,
\item $(x,yy') = ([\Res](x), y\otimes y'),$ for $x,y, y'\in K_0(R)$,
\item $(x x', y) = (x \otimes x', [\Res](y)),$ for $x, x', y\in K_0(R)$.
\end{enumerate}
\end{prop}
\begin{proof} Since $1=[P_{\emptyset}]$, where $\emptyset$ is the empty sequence,
and $P_{\emptyset}=\Bbbk$ as a module over $R(\emptyset)=\Bbbk$, the
first statement follows. When $i\not=j$, vectors $[P_i]$ and $[P_j]$ lie in mutually
orthogonal subspaces $K_0(R(i))$ and $K_0(R(j))$, so that $([P_i],[P_j])=0$. Also,
$$([P_i], [P_i]) =
\gdim({}_i R(i)_i) = \gdim(\Bbbk[x]) = (1-q^2)^{-1}.$$
Let $X\in R(\nu+\nu')\pmod$, $Y\in R(\nu)\pmod, $ and $Y'\in R(\nu')\pmod$. Then
\begin{eqnarray*}
& & ([X], [Y][Y'])= ([X], [\Ind_{\nu,\nu'} Y\otimes Y']) \\
& &  = \gdim (X^{\psi}\otimes_{R(\nu+\nu')}( {}_{\nu+\nu'}R_{\nu,\nu'} )
\otimes_{R(\nu)\otimes R(\nu')} Y\otimes Y' )  \\
& & = \gdim  (X^{\psi}(1_{\nu}\otimes 1_{\nu'})\otimes_{R(\nu)\otimes R(\nu')} Y\otimes Y' )
= ([\Res_{\nu,\nu'} X], [Y]\otimes [Y']),
\end{eqnarray*}
and statement 3 follows. A similar computation establishes 4.
\end{proof}

We next recall $\Q(q)$-algebras $'\mathbf{f}$ and $\mathbf{f}$
from~\cite[Section 1]{Lus4} and $_{\mc{A}}\mathbf{f}$, the integral
form of $\mathbf{f}$ (our $q$ is Lusztig's $v^{-1}$).
Algebra $'\mathbf{f}$ is a free associative
$\N[I]$-graded algebra on generators $\theta_i$. The degree of $\theta_i$ is $i$.
The tensor product $'\mathbf{f}\otimes '\mathbf{f}$ is equipped with a
algebra structure using the rule~\eqref{eq_quasi_commute}, and with a coalgebra
structure $r:\  '\mathbf{f}\lra '\mathbf{f}\otimes '\mathbf{f}$, determined by
the conditions $r(\theta_i)=\theta_i\otimes 1 + 1\otimes \theta_i$ and $r$
being an algebra homomorphism.

$'\mathbf{f}$ comes equipped with a bilinear form $(,)$ uniquely determined by the
same conditions as the ones in Proposition~\ref{prop-pairing-prop},
with comultiplication $r$ taking the place of comultiplication $[\Res]$ in
$K_0(R)$ and $\theta_i$, $\theta_j$ taking place of $[P_i]$, $[P_j]$.
It has the weight space decomposition
$$'\mathbf{f}= \bigoplus_{\nu\in \N[I]} {}'\mathbf{f}_{\nu}.$$
Let $\mc{I}$ be the radical of the bilinear form $(,)$. It is a two-sided
ideal of $'\mathbf{f}$ and one forms the quotient algebra $\mathbf{f}=
{}'\mathbf{f}/\mc{I}$, which also has the weight decomposition
$$\mathbf{f}= \bigoplus_{\nu\in \N[I]} \mathbf{f}_{\nu}.$$
The bilinear form and the comultiplication $r$ descend to the quotient algebra.
$'\mathbf{f}$ and $\mathbf{f}$ come with a $\Q(q)$-antilinear
involution $\bar{\quad}$ that takes $q^n$ to $q^{-n}$ and $\theta_i$ to $\theta_i$.

It is not hard to check that the elements
$$\theta_i \theta_j - \theta_i \theta_j \ \ {\rm for} \ \  i\cdot j=0 $$
 and
$$(q+q^{-1})\theta_i \theta_j \theta_i- \theta_i^2 \theta_j - \theta_j \theta_i^2 \ \
{\rm for} \ \ i\cdot j = -1$$
belong to the ideal $\mc{I}$. The quantum version of the
Gabber-Kac theorem says that $\mc{I}$ is generated by these elements over
all pairs of vertices $i\not= j$ of the graph $\Gamma$ (for instance, see
Theorem~33.1.3 in~\cite{Lus4}).

Define $\Af$ as the $\Z[q,q^{-1}]$-subalgebra generated by the divided
powers $\theta_i^{(a)}$, $i\in I,$ $a\in \N$.

\begin{prop} There is an injective homomorphism of $\Z[q,q^{-1}]$-algebras
$\gamma:\Af \lra K_0(R)$ that takes
$\theta_{i_1}^{(a_1)}\dots \theta_{i_k}^{(a_k)}$ to
$[P_{\ii}]$, where $\ii= i_1^{(a_1)}\dots i_k^{(a_k)}$.
This homomorphism converts the comultiplication $r$ of
$_{\mc{A}}\mathbf{f}$ into the comultiplication $[\Res]$ in  $K_0(R).$
It takes the bilinear form on $\Af$ to the bilinear form on $K_0(R)$:
$$(x,y) = (\gamma(x), \gamma(y)).$$
The bar-involution of $\Af$ goes to the bar-involution of $K_0(R)$ under $\gamma$.
\end{prop}
\begin{proof}
Start with the homomorphism of $\Q(q)$-algebras
$'\mathbf{f} \lra K_0(R)_{\Q(q)},$ where
$$K_0(R)_{\Q(q)}\define K_0(R)\otimes_{\Z[q,q^{-1}]} \Q(q),$$
defined by the condition that it takes $\theta_i$ to $[P_i].$
There are equalities in $K_0(R)_{\Q(q)}$:
\begin{eqnarray*}
 & &   [P_{ij}] = [P_{ji}] ,  \  \    i\cdot j=0 ,     \\
 & &   [P_{iji}] = [P_{i^{(2)}j}]+[P_{ji^{(2)}}], \  \   i\cdot j = -1,
\end{eqnarray*}
that come from isomorphisms of left
projective modules in Proposition~\ref{prop-iso-leftright}.
These equalities match the generators of the ideal $\mc{I}.$ Therefore,
the above homomorphism descends to a homomorphism
$$\gamma_{\Q(q)} \ : \ \mathbf{f} \lra K_0(R)_{\Q(q)} .$$
Under this homomorphism induction of projective $R$-modules
corresponds to the multiplication in $\Af$, so that
$$\gamma_{\Q(q)}(\theta_{i_1}\dots \theta_{i_k})=  [P_{i_1\dots i_k}].$$
Passing to the divided powers shows that
$$\gamma_{\Q(q)}(\theta_{i_1}^{(a_1)}\dots \theta_{i_k}^{(a_k)}) =
[P_{i_1^{(a_1)}\dots i_k^{(a_k)}}].$$
The bilinear forms on $\mathbf{f}$ and $K_0(R)_{\Q(q)}$ satisfy
the same properties, listed earlier, and these properties uniquely determine the
form on $\mathbf{f}$. Therefore, the homomorphism $\gamma_{\Q(q)}$
respects the bilinear forms:
$$(\gamma_{\Q(q)}(x), \gamma_{\Q(q)}(y)) = (x,y).$$
Since the bilinear form on $\mathbf{f}$ is non-degenerate, homomorphism
$\gamma_{\Q(q)}$ is injective. The bar-involution on $\mathbf{f}$ is
 $q$-antilinear and fixes each product element  $\theta_{i_1}\dots \theta_{i_k}$.
The bar-involution on $K_0(R)_{\Q(q)}$ is $q$-antilinear and fixes
$[P_{\ii}]$ for each $\ii.$ Therefore, $\gamma_{\Q(q)}(\overline{x}) =
 \overline{\gamma_{\Q(q)}(x)}$ for all $x\in \mathbf{f}$.

The image of the restriction of $\gamma_{\Q(q)}$ to $\Af$ lies in
$K_0(R)$, therefore we get a homomorphism $\gamma: \Af \lra K_0(R)$
by restriction. This homomorphism is injective and satisfies all the properties
stated in the proposition.
\end{proof}

We will prove in the next section that $\gamma$ is an isomorphism.

%
\subsection{Surjectivity of $\gamma$}
\label{subsec-injectivity}
%

In this section we closely follow~\cite[Chapter 5]{KleBook}; all results
there transfer directly to our case.

For $M$ in $R(\nu)\dmod$ and $i\in I$ let
$$\Delta_{i}M= (1_{\nu-i}\otimes 1_i)M = {}_{\nu-i,i}R_{\nu}\otimes_{R(\nu)} M,$$
and, more generally,
$$\Delta_{i^n}M =  (1_{\nu-ni}\otimes 1_{ni})M = {}_{\nu-ni,ni}R_{\nu}
 \otimes_{R(\nu)} M.$$
We view $\Delta_{i^n}$ as the functor into the category $R(\nu-ni)\otimes R(ni)\dmod$.
There are functorial isomorphisms
\begin{equation} \label{eq_funs}
\HOM_{R(\nu)} (\Ind_{\nu-ni,ni}N\otimes L(i^n), M) \cong
\HOM_{R(\nu-ni)\otimes R(ni)}(N\otimes L(i^n), \Delta_{i^n}M),
\end{equation}
for $M$ as above and $N\in R(\nu-ni)\dmod$. The following lemma is obvious.

\begin{lem} \label{lemma_shuffle2}
$$\chr (\Delta_{i^n}M)= \sum_{\jj\in \seq(\nu-ni)}   \chr(M, \jj i^n)\cdot \jj, $$
where we view $\Delta_{i^n}M$ as a module over the subalgebra $R(\nu-ni)$
of $R(\nu-ni)\otimes R(ni)$.
\end{lem}

Let $\epsilon_i(M) = {\rm max}\{n\ge 0 | \Delta_{i^n}M\not= 0\}$.
This number is the length of the longest tail of $i$'s in sequences $\kk$ with
$1_{\kk}M \not= 0$.

\begin{lem} If $M\in R(\nu)\dmod$ is irreducible, and $N\otimes L(i^n)\{r\}$
is an irreducible submodule of $\Delta_{i^n}(M)$ for some $0\le n \le \epsilon_i(M)$
and $r\in \Z$, then $\epsilon_i(N)  =\epsilon_i(M) - n.$
\end{lem}
\begin{proof} This is our analogue of Lemma~5.1.2 of~\cite{KleBook} and the proof
is essentially the same. Let $\epsilon=\epsilon_i(M)$.
Clearly, $\epsilon_i(N) \le \epsilon_i(M)-n$.
Isomorphisms~\eqref{eq_funs} and the irreducibility of $M$ imply that it is a
quotient of $\Ind_{\nu-ni, ni} N \otimes L(i^n)\{r\}$. By exactness of $\Delta_{i^{\epsilon}}M$,
we get that $\Delta_{i^{\epsilon}}(M)\not= 0$ is a quotient of
$$ \Delta_{i^{\epsilon}}(\Ind_{\nu-ni, ni}N\otimes L(i^n))\{r\}.$$
Hence, the latter module is nonzero, and the inequality $\epsilon_i(N)\ge \epsilon_i(M)-n$
follows from the Shuffle lemma~\ref{lemma_shuffle}.
\end{proof}

\begin{lem} Suppose $N\in R(\nu)\dmod$ is irreducible and $\epsilon_i(N)=0$.
Let $M=\Ind_{\nu, ni} N\otimes L(i^n)$. Then
\begin{enumerate}
\item $\Delta_{i^n}M\cong N\otimes L(i^n)$,
\item ${\rm hd} M$ is irreducible and $\epsilon_i({\rm hd} M)= n$,
\item all other composition factors $L$ of $M$ have $\epsilon_i(L)<n$.
\end{enumerate}
\end{lem}
\begin{proof} This is the analogue of Lemma~5.1.3 in~\cite{KleBook} for algebras $R(\nu)$.

1) is immediate from the Shuffle lemma and Lemma~\ref{lemma_shuffle2}.

2) From~\eqref{eq_funs} we see that a copy of $N\otimes L(i^n)$, possibly with
a grading shift, appears in
$\Delta_{i^n}Q$ for any non-zero quotient $Q$ of $M$, including direct summands
of ${\rm hd}M$. Part 1, however, implies that $N\otimes L(i^n)$ appears only once
in $\Delta_{i^n}M$, so that ${\rm hd}M$ is irreducible.

3) From part 2) we have $\Delta_{i^n}(M)= \Delta_{i^n}({\rm hd}M)$, so that
$\Delta_{i^n}(L)=0$ for any other composition factor of $M$, since $\Delta_{i^n}$
is exact.
\end{proof}

\begin{lem} Let $M\in R(\nu)\dmod$ be irreducible and $\epsilon=\epsilon_i(M)$.
Then $\Delta_{i^{\epsilon}}M$ is isomorphic to $N\otimes L(i^{\epsilon})$
for some irreducible $N\in R(\nu-\epsilon i)\dmod$ with $\epsilon_i(N)=0$.
\end{lem}

\begin{proof} The proof is identical to that of Lemma~5.1.4 in~\cite{KleBook}.
\end{proof}

\begin{lem} \label{lemma_515}
Let $N\in R(\nu)\dmod$ be irreducible and $M=\Ind_{\nu,ni}N\otimes L(i^n)$.
Then ${\rm hd}M$ is irreducible, $\epsilon_i({\rm hd}M)= \epsilon_i(N)+n$, and
all other composition factors $L$ of $M$ have $\epsilon_i(L)<\epsilon_i(N)+n$.
\end{lem}

\begin{proof} Same as proof of Lemma~5.1.5 in~\cite{KleBook}.
\end{proof}

\begin{prop} For any irreducible $M\in R(\nu)\dmod$ and $0\le n \le \epsilon_i(M)$,
${\rm soc} \Delta_{i^n}M$ is an irreducible $R(\nu-ni)\otimes R(ni)$-module of
the form $L\otimes L(i^n)$ with $\epsilon_i(L)=\epsilon_i(M)-n$.
\end{prop}

\begin{proof} Same as proof of Theorem~5.1.6 in~\cite{KleBook}.
The analogue of the Kato theorem in our framework is stated below
(this theorem appears in the proof of Theorem~5.1.6).
\end{proof}

\begin{prop} Let $\mu$ be a composition of $n$.
\begin{enumerate}
\item The module $L(i^n)$ over the nilHecke algebra $R(ni)$ is the only graded
irreducible module, up to isomorphism and graded shifts.
\item All composition factors of $\Res^n_{\mu}L(i^n)$ are isomorphic
to $L(i^{\mu_1})\otimes \dots \otimes L(i^{\mu_r})$, up to grading shifts, and
${\rm soc}(\Res^n_{\mu} L(i^n))$ is irreducible.
\item ${\rm soc}(\Res^n_{n-1} L(i^n))\cong L(i^{n-1})$, up to a grading shift.
\end{enumerate}
\end{prop}
Here $\Res^n_{\mu}$ denotes the restriction to the parabolic nilHecke subalgebra
$\BNC_{\mu}\cong \BNC_{\mu_1}\otimes \dots \otimes \BNC_{\mu_r}$.
The proof in~\cite{KleBook} works in this case as well, with the equivalent
of Lemma~4.3.1 being Lemma~\ref{lem-prekato}.  $\square$

Let $e_i = \Res^{\nu-i, i}_{\nu-i}\circ \Delta_i$ be the functor of composition of $\Delta_i$
with the restriction from $R(\nu-i)\otimes R(i)$ to $R(\nu-i)$. Then
$\epsilon_i(M)= {\rm max}\{n\ge0| e^n_i M \not=0\}$ and
$$ \Res^{\nu}_{\nu-i} M = \bigoplus_{i\in I} e_i M.$$

\begin{cor} \label{cor_517}
Let $M\in R(\nu)\dmod$ be irreducible with $\epsilon_i(M)>0$.
Then ${\rm soc}(e_i M)$ is irreducible and $\epsilon_i({\rm soc}(e_i M))=\epsilon_i(M)-1$.
Socles of $e_iM$ are pairwise non-isomorphic for different $i\in I$.
\end{cor}

Proof is the same as for Corollaries~5.1.7 and 5.1.8 in~\cite{KleBook}. $\square$

For an irreducible $M\in R(\nu)\dmod$ define
\begin{equation}
 \tilde{e}_iM := {\rm soc} (e_i M), \ \
  \tilde{f}_i M := \textrm{hd ind}^{\nu+i}_{\nu,i} M\otimes L(i).
\end{equation}
The module $\tilde{f}_iM$ is irreducible by Lemma~\ref{lemma_515}, while
$\tilde{e}_iM$ is irreducible or $0$ by Corollary~\ref{cor_517}, and
\begin{equation*}
 \epsilon_i(M) = {\rm max}\{ n\ge 0| \tilde{e}_i^n M \not= 0 \},
\ \  \epsilon_i(\tilde{f}_i M) = \epsilon_i(M)+1 .
\end{equation*}

In the statements below, isomorphisms of simple modules are allowed to
be homogeneous (not necessarily degree-preserving).

\begin{lem} For an irreducible $M\in R(\nu)\dmod$ we have
\begin{eqnarray}
 {\rm soc} \Delta_{i^n}M & \cong & (\tilde{e}^n_iM)\otimes L(i^n), \\
  \textrm{hd ind}_{\nu,ni} (M\otimes L(i^n)) & \cong & \tilde{f}^n_i M.
\end{eqnarray}
\end{lem}

\begin{lem} For an irreducible $M\in R(\nu)\dmod$ the socle of $e_i^n M$ is
isomorphic to $\tilde{e}^n_iM^{\oplus [n]!}\{-\frac{n(n-1)}{2}\}$.
\end{lem}

The proofs are equivalent to those of Lemmas~5.2.1 and 5.2.2 in~\cite{KleBook}.
$\square$

\begin{lem} For irreducible modules $M\in R(\nu)\dmod$ and $N\in R(\nu+i)\dmod$
we have $\tilde{f}_i M \cong N$ if and only if $\tilde{e}_i N\cong M$.
\end{lem}

Proof follows that of Lemma~5.2.3 in~\cite{KleBook}. $\square$

\begin{cor} Let $M,N\in R(\nu)\dmod$ be irreducible. Then
$\tilde{f}_iM \cong \tilde{f}_iN$ if and only if $M\cong N$.
Assuming $\epsilon_i(M), \epsilon_i(N)>0$, $\tilde{e}_i M \cong \tilde{e}_i N$
if and only if $M\cong N$.
\end{cor}

The character $\chr(M)$ of a finite-dimensional representation $M\in R(\nu)\dmod$
takes values in $\Z[q,q^{-1}]\seq(\nu)$, the free $\Z[q,q^{-1}]$-module generated by
$\seq(\nu)$, and descends to a homomorphism from the Grothendieck group
$G_0(R(\nu))$ to $\Z[q,q^{-1}]\seq(\nu)$.

\begin{thm} The character map
$$ \chr: G_0(R(\nu)) \lra \Z[q,q^{-1}]\seq(\nu) $$
is injective.
\end{thm}

Equivalently, the characters of irreducible modules (one from each
equivalence class up to grading shifts) are linearly independent functions
on $\seq(\nu)$. The proof is identical to that of Theorem~5.3.1 in~\cite{KleBook}.
Note that in our case the character of a finite-dimensional graded module is a
function on sequences with values in $\Z[q,q^{-1}]$, while in the non-graded
case of~\cite{KleBook} its a function on sequences taking values in $\Z$.
This discrepancy has no effect on the proof. $\square$

Passing to the fraction field $\Q(q)$ of $\Z[q,q^{-1}]$ and dualizing the
map $\chr$, which then becomes the composition
$$\Q(q)\seq(\nu) \lra \mathbf{f}_{\nu} \stackrel{\gamma_{\Q(q)}}{\lra}
K_0(R(\nu))_{\Q(q)},$$
we conclude that $\gamma_{\Q(q)},$ restricted to weight $\nu$,
is a surjective map of $\Q(q)$-vector spaces.
We have already  observed that $\gamma$ and $\gamma_{\Q(q)}$ are
injective. By summing over all weights, we obtain the following result.

\begin{prop} $\gamma_{\Q(q)}: \mathbf{f}\lra K_{0}(R)_{\Q(q)}$ is
an isomorphism.
\end{prop}

Therefore, the number of isomorphism classes of (graded) simple $R(\nu)$-modules is the
same for any field $\Bbbk$.

\begin{cor} A (graded) irreducible $R(\nu)$-module is absolutely irreducible,
for any $\Gamma$, $\Bbbk$ and weight $\nu$.
\end{cor}

Next, assume that $\Gamma$ is finite. Choose a total order on $I$,
$i(0)<i(1)\dots <i(k-1), $  $k=|I|$. For $r>k$ define $i(r)=i(r')$ where $r'$ is
the residue of $r$ modulo $k$.
Fix one representative $S_b$ from each isomorphism class $b$ of irreducible $R(\nu)$-modules,
up to grading shifts. Recall that we denoted this set of isomorphism classes by
$\mathbf{B}'_{\nu}$.  For all $\nu$, to each $b\in \mathbf{B}'_{\nu}$ assign the
following sequence $Y_b=y_0y_1\dots$ of nonnegative integers:
$y_0=\epsilon_{i(0)}(M)$, and  let $M_1=\tilde{e}^{y_0}_{i(0)} M$.
Inductively,
$y_r=\epsilon_{i(r)}(M_r)$, and  $M_{r+1}=\tilde{e}^{y_r}_{i(r)} M_r$.
Note that $y_0+y_1+\dots = |\nu|$ and only finitely many terms in the
sequence are non-zero. Introduce a lexicographic order on sequences of non-negative
integers: $y_0y_1 \dots > z_0z_1 \dots$ if, for some $t$, $y_0=z_0,$ $y_1=z_1$, \dots
$y_{t-1}=z_{t-1}$ and $y_t > z_t.$ This order induces a total order on
$\mathbf{B}'_{\nu}$, by $b>c$ iff $Y_b> Y_c.$
To each sequence $Y_b=y_0y_1\dots $ we assign the projective $R(\nu)$-module
$P_{Y^r_b}$ associated to the divided powers sequence
$Y^r_b=\dots i(2)^{(y_2)}i(1)^{(y_1)}i(0)^{(y_0)}$ (the order of $y$'s is reversed).

\begin{prop} $\HOM(P(Y_b), S_c)=0$ if $b<c$ and
$\HOM(P(Y_b), S_b) = \Bbbk$.
\end{prop}

This follows from the previous results and implies that the image $[P]$
of any (graded) projective $R(\nu)$-module in the Grothendieck group $K_0(R(\nu))$
can be written as a linear combination, with coefficients in $\Z[q,q^{-1}]$,
 of images of divided powers projectives
$[P_{\theta}]$, for divided power sequences $\theta$ of the form $Y^r_b$.
Therefore, $\gamma: \Af\lra K_0(R(\nu))$ is surjective. Since, $\gamma$ is
also injective, it is an isomorphism. The case of an infinite $\Gamma$ follows
by taking the direct limit of its finite subgraphs.
This concludes the proof of the Theorem~\ref{thm_first}
stated in the introduction.

It would be interesting to find out if the theorem remains valid for rings $R(\nu)$
over $\Z$ rather than over a field $\Bbbk$.

\vspace{0.1in}

For each divided power $i^{(a)}$
we have the corresponding projective $P_{i^{(a)}}$.
Induction with this projective is an exact functor, denoted
$\mathcal{F}_i^{(a)}$, from $R(\nu)\dmod$ to $R(\nu+ai)\dmod$.
Summing over all $\nu$, form the functor
$$ \mathcal{F}_i^{(a)} \  :  \   R\dmod \lra R\dmod. $$
This functor restricts to the subcategory $R\pmod$ of the category of projective modules.
To any divided power sequence $\theta=i_1^{(a_1)}\dots i_r^{(a_r)}$
associate the functor
$$\mathcal{F}_{\theta} = \mathcal{F}_{i_1}^{(a_1)}\circ \dots
\circ\mathcal{F}_{i_r}^{(a_r)}$$
on $R\dmod$. To a finite sum $\sum_k u_k {\theta(k)}$ where
$u_k \in \N[q,q^{-1}]$ and $\theta(k)$ are divided powers sequences associate
the direct sum of shifted copies of $\mathcal{F}_{\theta(k)}$:
$$ \bigoplus_{k} \mathcal{F}_{\theta(k)}^{\oplus u_k} .$$

\begin{thm} For any relation
$$\sum_k u_k \theta(k)= \sum_{\ell} v_{\ell} \theta'(\ell)$$
in $\Af$ with positive coefficients $u_k, v_{\ell} \in \N[q,q^{-1}]$ there
is an isomorphism of projectives
$$ \bigoplus_{k} P_{\theta(k)}^{\oplus u_k} \cong
 \bigoplus_{\ell} P_{\theta'(\ell)}^{\oplus v_{\ell}}$$
inducing an isomorphism of functors
$$  \bigoplus_{k} \mc{F}_{\theta(k)}^{\oplus u_k} \cong
\bigoplus_{\ell} \mc{F}_{\theta'(\ell)}^{\oplus v_{\ell}}.$$
\end{thm}

This result follows immediately from the earlier ones. $\square$

We conclude that any relation in $\Af$ lifts to an isomorphism of functors. It is
natural to view the category $R\pmod$, as well as the category of induction functors
on $R\dmod$ it gives rise to, as a categorification of $\Af$, the integral form
of the quantum universal enveloping algebra of the negative half of the simply-laced
Kac-Moody algebra associated to the graph $\Gamma$.

\vspace{0.1in}

The semilinear "hom" form $(,)'$ on $K_0(R)$ defined by 
$$([P],[Q])' := \gdim \HOM (P,Q)$$
is related to the "tensor product" bilinear form $(,)$ given by~\eqref{eq_bil_pair2} via 
$$ (x,y)' = (\overline{x},y).$$ 
Indeed, by surjectivity of $\gamma$, it suffices to check this relation for $x=[P_{\ii}]$ and 
$y=[P_{\jj}]$, in which case both sides are equal to the graded dimension of 
${}_{\ii} R(\nu)_{\jj}$. 

The involution $\sigma$ of $R(\nu)$ defined in Section~\ref{subsec-definitions}
induces a self-equivalence of $R(\nu)\dmod$ which takes projective
$P_{i_1^{(a_1)}\dots i_r^{(a_r)}}$ to $P_{i_r^{(a_r)}\dots i_1^{(a_1)}}$.
The induced map $[\sigma]$ on the Grothendieck group $K_0(R)$ coincides,
under the isomorphism $\gamma$, the $q$-linear anti-involution of
$\Af$ that fixes each $\theta_i^{(a)}$.

On the category $R(\nu)\fmod$ we have the contravariant duality functor
which takes a finite-dimensional module $M$ to its vector space dual $M^{\ast \psi}$
twisted by the antiinvolution $\psi$. This duality functor leaves invariant the
character evaluated at $q=1$:
$$ \chr(M^{\ast\psi})_{q=1} = \chr(M)_{q=1} .$$
Therefore, the contravariant duality preserves simples, up to overall shift:
$$ S_b^{\ast \psi} \cong S_b\{r\}.$$
Given $\ii\in \seq(\nu)$,
we have $\chr(S_b, \ii) \in \Z[q^2,q^{-2}]$ or $\chr(S_b, \ii) \in q\Z[q^2,q^{-2}]$
for parity reasons (more generally, this is true for any indecomposable object
of $R(\nu)\dmod$ and can be used to decompose $R(\nu)\dmod$ into the
direct sum of two subcategories). Then  $\chr(S_b^{\ast\psi}, \ii) \in \Z[q^2,q^{-2}]$
if the same is true for $\chr(S_b, \ii)$, and  $\chr(S_b^{\ast\psi}, \ii) \in q\Z[q^2,q^{-2}]$
if  $\chr(S_b, \ii) \in q\Z[q^2,q^{-2}]$. Hence, the shift $r$ is an even number.

From now on we redefine $S_b$ by shifting its grading by $\frac{r}{2}$.
We have $ S_b^{\ast \psi} \cong S_b$ as graded modules. This normalization of $S_b$
does not depend on the choice of $\ii$. The character of $S_b$ is bar-invariant:
$$ \overline{\chr(S_b, \ii)} = \chr(S_b, \ii)$$
for all $\ii \in \seq(\nu)$, where $\overline{q}=q^{-1}$.  Extending the bar-involution
to $\Z[q,q^{-1}]\seq(\nu)$ by $\overline{\ii}=\ii$, we have
$\overline{\chr(S_b)}= \chr(S_b)$.

This canonical (balanced) choice of grading for $S_b$ allows us to fix
the grading on indecomposable projective $P_b$ so that the quotient map
$P_b\lra S_b$ is grading-preserving. In this way we obtain a basis $\{[P_b]\}$  
in $\Af$ which depends only on the characteristic of $\Bbbk$. Both the 
multiplication and the comultiplication in this basis have coefficients in 
$\N[q,q^{-1}]$. An example 
below shows this basis to be different from the Lusztig-Kashiwara basis 
when $\Gamma$ is an odd length cycle and $\Bbbk$ has any characteristic, 
and when $\Gamma$ is a  cycle and $\Bbbk$ has characteristic $2$.

%
\subsection{Tight monomials and indecomposable projectives}
%

Following Lusztig~\cite{Lus5}, we say that a monomial
$\theta=\theta_1^{(a_1)}\dots \theta_k^{(a_k)}$ is \emph{tight} if
it belongs to the canonical basis $\mathbf{B}$ of $\Af$.
It follows from the properties of the canonical basis that a monomial
$\theta$ is tight if and only if $(\theta, \theta) -1 \in  q \N[q]$
(or see~\cite[Proposition 3.1]{Reineke}).

\begin{prop} If a monomial $\theta$ is tight, the projective module $P_{\theta}$ is
indecomposable.
\end{prop}
\begin{proof} Tightness of $\theta$ implies that $\HOM(P_{\theta},P_{\theta})$
is a $\Z_+$-graded $\Bbbk$-vector space which is one-dimensional in degree $0$. Therefore,
any degree $0$ endomorphism of $P_{\theta}$ is a multiple of the identity,
and $P_{\theta}$ is indecomposable. This argument works even over $\Z$.
\end{proof}

\begin{example} When the graph $\Gamma$ consists of a single vertex $i$,
the weight space $\Af_{mi}$ is a rank one free $\Z[q,q^{-1}]$-module generated
by $\theta_i^{(m)}$. The map $\gamma$ takes it to $[P_{i^{(m)}}]$, the
generator of the Grothendieck group $K_0(R(mi))$. Projective module $P_{i^{(m)}}$
is indecomposable.
\end{example}

\begin{example} Let $\Gamma=\xy
  (-5,0)*{\circ};
  (5,0)*{\circ}
  **\dir{-};
  (-5,2.5)*{\scs i};
  (5,2.5)*{\scs j};
  \endxy$. Tight monomials
$\theta_i^{(a)}\theta_j^{(b)}\theta_i^{(c)}$ ($a,b,c\in\N$, $b\ge a+c$) and
$\theta_j^{(c)}\theta_i^{(b)}\theta_j^{(a)}$ ($a,b,c\in\N$, $b\ge a+c$),
with the identification
$$\theta_i^{(a)}\theta_j^{(a+c)}\theta_i^{(c)}=
  \theta_j^{(c)}\theta_i^{(a+c)}\theta_j^{(a)}, $$
constitute the canonical basis $\mathbf{B}$ of $\Af$,
see~\cite[Example 14.5.4]{Lus4}. Therefore, images of indecomposable projectives
$P_{i^{(a)}j^{(b)}i^{(c)}}$, $b\ge a+c$ and $P_{j^{(c)}i^{(b)}j^{(a)}}$,
$b> a+c$, constitute a basis in the free $\Z[q,q^{-1}]$-module $K_0(R)$.
Any indecomposable projective in $R\dmod$ is isomorphic to one of the above,
up to a grading shift.
Indecomposables $P_{i^{(a)}j^{(a+c)}i^{(c)}}$ and $P_{j^{(c)}i^{(a+c)}j^{(a)}}$ are
isomorphic.
\end{example}

\begin{example} 
Let $\Gamma = \xy 0;/r.12pc/:  (-10,0)*{\circ }="1";
 (5,10)*{\circ }="3"; (-5,10)*{\circ }="2";
  (5,-10)*{\circ }="5";
  (-5,-10)*{\circ }="6";
 (10,0)*{\circ }="4";
 "1";"2" **\dir{-};
 "2";"3" **\dir{-};
 "3";"4" **\dir{-};
 "4";"5" **\dir{-};
 "5";"6" **\dir{-};
 "6";"1" **\dir{-};
\endxy$ be a cycle with $n \geq 3$ vertices.  Label the vertices clockwise by $1,2, \dots, n$ and
let $\ii=12\dots n$.  Then $\HOM(P_{\ii\ii},P_{\ii\ii})$ is $\Z_+$-graded and
$\Hom(P_{\ii\ii}, P_{\ii\ii})$ is 2-dimensional with the basis 
$\{ {}_{\ii\ii}1_{\ii\ii},\alpha\}$, where 
\[
 \alpha =
 \xy
  (0,0)*{\includegraphics[scale=0.5]{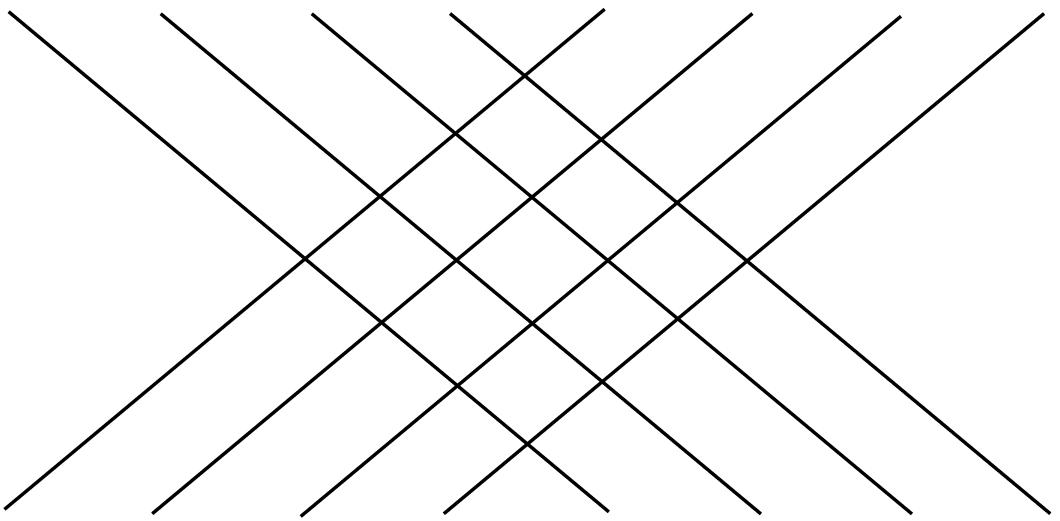}};
  (-27,-15)*{\scs 1};(-19,-15)*{\scs 2};(-12,-15)*{\scs 3};(-8,-15)*{\scs \cdots};(-4,-15)*{\scs n};
  (4,-15)*{\scs 1};(12,-15)*{\scs 2};(19,-15)*{\scs 3};(23,-15)*{\scs \cdots};(27,-15)*{\scs n};
 \endxy
\]
A computation shows that $\deg(\alpha)=0$ and
\begin{equation*}
 \alpha^2 = \left\{ \begin{array}{cl}
                      0 & \text{if $n$ is odd,} \\
                      -2 \alpha & \text{if $n$ is even;}
                    \end{array}
 \right.
\end{equation*}
implying that ${_{\ii\ii}1_{\ii\ii}}$ is the only idempotent in $P_{\ii\ii}$ if
$n$ is odd, or if $n$ is even and ${\rm char} (\Bbbk) =2$. Under these assumptions,
$P_{\ii\ii}$ is indecomposable, but
\begin{equation*}
 \big([P_{\ii\ii}],[P_{\ii\ii}]\big) \in 2 + q\N[q] \neq 1+q \N[q], 
\end{equation*}
and $[P_{\ii\ii}]$ is not a canonical basis element. For $\Gamma$ an odd length
cycle we found an indecomposable projective $P_{\ii\ii}$ whose image in the
Grothendieck group is not a canonical basis vector, while being invariant under the bar 
involution: $\overline{[P_{\ii\ii}]}= [P_{\ii\ii}]$.

When $n$ is even and ${\rm char} (\Bbbk) \neq 2$, $\alpha_0=-\frac{\alpha}{2}$ is
an idempotent in $\Hom(P_{\ii\ii},P_{\ii\ii})$, and
$P_{\ii\ii} \cong P_{\ii\ii}\alpha_0 \oplus P_{\ii\ii}(1-\alpha_0)$
is isomorphic to the direct sum of two indecomposable projectives. 
Furthermore, $\alpha = \beta_1 \beta_0$ where
\[
 \beta_1 =
 \xy
  (0,0)*{\includegraphics[scale=0.5]{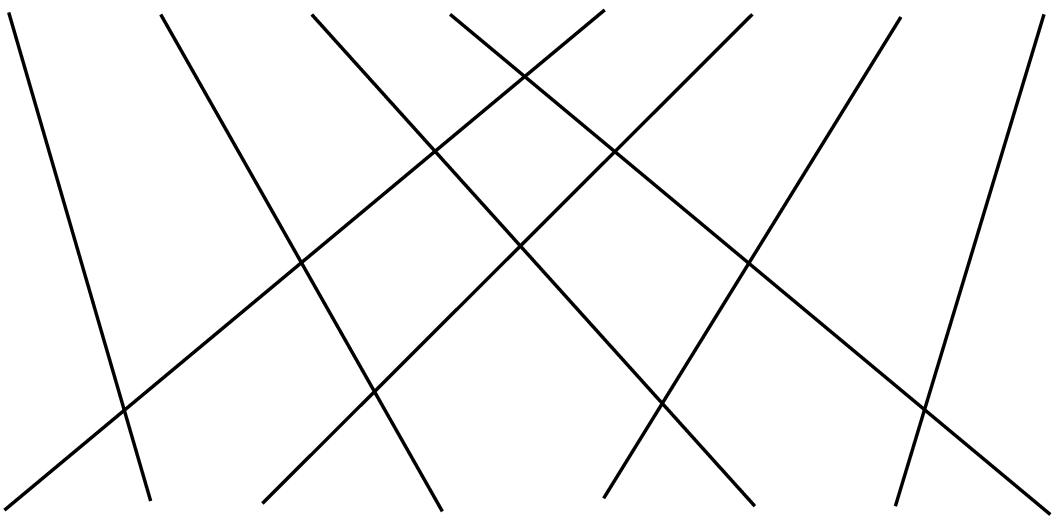}};
  (-27,-15)*{\scs 1};(-19,-15)*{\scs 1};(-13,-15)*{\scs 2};(-4,-15)*{\scs 2};
  (4,-15)*{\scs 3};(12,-15)*{\scs 3};(19,-15)*{\scs n};(15.5,-15)*{\scs \cdots};(27,-15)*{\scs n};
  (-23.5,-10)*{\bullet};(-11,-10)*{\bullet};(5.3,-10)*{\bullet};(19.5,-10)*{\bullet};
 \endxy
\]
\[
 \beta_0 =
 \xy
  (0,0)*{\includegraphics[scale=0.5]{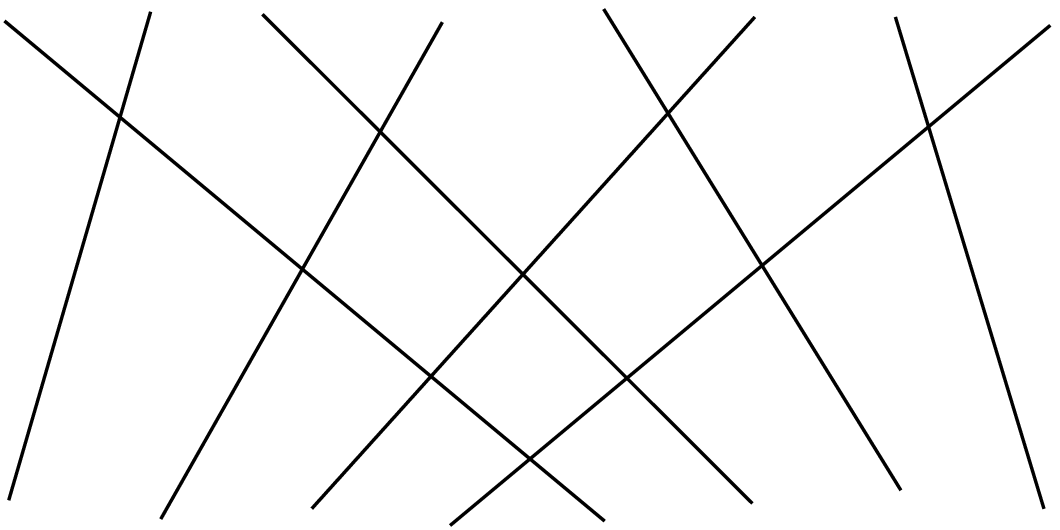}};
 (-27,-15)*{\scs 1};(-19,-15)*{\scs 2};(-12,-15)*{\scs 3};(-8,-15)*{\scs \cdots};(-4,-15)*{\scs n};
  (4,-15)*{\scs 1};(12,-15)*{\scs 2};(19,-15)*{\scs 3};(23,-15)*{\scs \cdots};(27,-15)*{\scs n};
 \endxy
\]
Module homomorphisms 
\begin{eqnarray*}
  \beta_1 &\maps & P_{\ii\ii}  \longrightarrow P_{1^{(2)}2^{(2)}\dots n^{(2)}}, \\
  \beta_0 &\maps & P_{1^{(2)}2^{(2)}\dots n^{(2)}}  \longrightarrow P_{\ii\ii},
\end{eqnarray*}
induced by these elements via right multiplication have degree $0$,  
and $\beta_0\beta_1 = - 2 \cdot \Id$, so that $P_{\ii\ii}\alpha_0 \cong 
P_{1^{(2)}2^{(2)}\dots n^{(2)}}$. 
\end{example}

%
\subsection{A conjecture on categorification of irreducible representations}
\label{subsec-conj-irr}
%

Choose $\lambda\in \N[I],$ $\lambda=\sum \lambda_i \cdot i ,$ $i\in I$.
Let $R(\nu;\lambda)$ be the quotient ring of $R(\nu)$ by the ideal
generated by all diagrams of the form
\[
\xy (0,18)*{}; (0,0)*{\up{i_1}}; (5,0)*{\up{\;\; i_2}};
  (10,0)*{\up{ }}; (15,0)*{\dots}; (20,0)*{\up{\;\;i_m}};
  (0,0)*{\bullet}+(-4,0)*{\lambda_{i_1}};
 \endxy
\]
where $i_1\dots i_m\in \seq(\nu)$ and the leftmost string has
$\lambda_{i_1}$ dots on it.  The ring $R(\nu; \lambda)$ inherits a
grading from $R(\nu)$. These quotient rings should be the analogues
of the Ariki-Koike cyclotomic Hecke algebras in our framework.  Let
$$ R(\ast;\lambda) \define \bigoplus_{\nu\in \N[I]} R(\nu;\lambda),$$
and switch from $\Z$ to a field $\Bbbk$.
We expect that, for sufficiently nice $\Gamma$ and $\Bbbk$,
the category of graded modules over $R(\ast;\lambda)$
categorifies the integrable irreducible $U_q(\mf{g})$-representation
$V_{\lambda}$ with the highest weight $\lambda$.
 Let $R(\nu;\lambda)\pmod$ be the
category of finitely-generated graded projective left $R(\nu;\lambda)$-modules and
$$ R(\ast;\lambda)\pmod \define \bigoplus_{\nu\in \N[I]} R(\nu;\lambda)\pmod .$$
There should exist an isomorphism
$$K_0(R(\ast;\lambda)) \cong V_{\Z,\lambda},$$
where $V_{\Z,\lambda}$ is an integral version of
$V_{\lambda}$, a free $\Z[q,q^{-1}]$-module spanned by the compositions
of divided differences $F_i^{(a)}$ applied to the highest weight vector $v_{\lambda} \in
V_{\lambda}$. Under this isomorphism indecomposable projectives should
correspond to canonical basis vectors in $V_{\lambda}$.
The action of $E_i^{(a)}$ and $F_i^{(a)}$ should lift to exact functors
$\mc{E}_i^{(a)}$ and $\mc{F}_i^{(a)}$ between categories
$R(\nu;\lambda)\pmod$ and $R(\nu+ai;\lambda)\pmod$ as well as the
categories $R(\nu;\lambda)\dmod$ and $R(\nu+ai;\lambda)\dmod$
of all finitely-generated graded modules. These functors
$\mc{E}_i^{(a)}$ and $\mc{F}_i^{(a)}$ will be direct summands of the induction
and restriction functors between $R(\nu;\lambda)$ and $R(\nu+ai;\lambda)$-modules,
defined \'a la Ariki. We expect them to be biadjoint, up to grading shifts.




\vspace{0.1in} 

\noindent 
M.K.: {\sl \small School of Mathematics, Institute for Advanced Study, Princeton, NJ 08540, and} 
{ \sl \small Department of Mathematics, Columbia University, New York, NY 10027}
\noindent 
  {\tt \small email: khovanov@math.columbia.edu}

\vspace{0.1in} 

\noindent 
A.L:  { \sl \small Department of Mathematics, Columbia University, New York, NY 10027}
\noindent 
  {\tt \small email: lauda@math.columbia.edu}

\end{document}